\documentclass{amsart}
\usepackage[all]{xy}
\usepackage{amsmath}
\usepackage{amsfonts,amssymb}
\usepackage{latexsym}
\RequirePackage[T1]{fontenc}

\numberwithin{equation}{section}

\newcommand{\al}{{\alpha}}
\newcommand{\la}{{\lambda}}
\newcommand{\io}{{\iota}}

\newcommand{\ga}{{\gamma}}

\newcommand{\ka}{{\kappa}}
\newcommand{\vr}{{\varrho}}
\newcommand{\de}{{\delta}}
\newcommand{\D}{{\Delta}}
\newcommand{\n}{{\nabla}}
\newcommand{\G}{{\Gamma}}

\newcommand{\Om}{{\Omega}}

\newcommand{\si}{{\sigma}}

\newcommand{\bz}{\bar{z}}

\newcommand{\bA}{\bar{A}}

\newcommand{\brd}{\bar{d}}

\newcommand{\inc}{\,{\rm inc}\,}

\newtheorem{dfn}{Definition}[section]
\newtheorem{prop}[dfn]{Proposition}
\newtheorem{cor}[dfn]{Corollary}

\newtheorem{thm}[dfn]{Theorem}
\newtheorem{lem}[dfn]{Lemma}

\newcommand{\Linf}{L_{\infty}}
\newcommand{\hotimes}{{\hat{\otimes}}}
\newcommand{\bul}{{\bullet}}
\newcommand{\brarrow}{\succ\rightarrow}
\newcommand{\blarrow}{\leftarrow\prec}
\newcommand{\bbrarrow}{\succ\succ\rightarrow}
\newcommand{\bblarrow}{\leftarrow\prec\prec}
\newcommand{\erarrow}{\stackrel{\sim}{\rightarrow}}

\newcommand{\bB}{{\bf B}}
\newcommand{\id}{\mathop{\rm{id}}}

\newcommand{\h}{\hbar}

\newcommand{\K}{\mathbb{K}}
\newcommand{\R}{\mathbb{R}}

\newcommand{\C}{\mathbb{C}}

\newcommand{\BV}{\,\Big|_{V}}
\newcommand{\BE}{\,\Big|_{E}}
\newcommand{\BTO}{\,\Big|_{T^{0,1}}}

\newcommand{\cL}{\mathcal{L}}
\newcommand{\cM}{\mathcal{M}}
\newcommand{\cS}{\mathcal{S}}
\newcommand{\cK}{\mathcal{K}}
\newcommand{\cR}{\mathcal{R}}
\newcommand{\cH}{\mathcal{H}}

\newcommand{\Tp}{\mathcal{T}}
\newcommand{\Dp}{\mathcal{D}}
\newcommand{\Ef}{\mathcal{A}}
\newcommand{\Jp}{\mathcal{J}}
\newcommand{\Bu}{\mathcal{B}}
\newcommand{\Op}{\mathcal UE}
\newcommand{\mO}{\mathcal O}

\newcommand{\tcK}{{\mathcal K}^{tw}}
\newcommand{\tcS}{{\mathcal S}^{tw}}

\newcommand{\SE}{\hat{S}(E^{\vee})}

\newcommand{\ETp}[1]{~\!\!^ET_{poly}^{#1}}
\newcommand{\ECp}[1]{~\!\!^E\!C^{poly}_{#1}}
\newcommand{\EDp}[1]{~\!\!^E\!D_{poly}^{#1}}
\newcommand{\EA}[1]{~\!\!^E\!A_{#1}}
\newcommand{\EJp}[1]{~\!\!^E\!\!J^{poly}_{#1}}
\newcommand{\Edif}{~\!\!^E\!d}
\newcommand{\Elie}{~\!\!^E\!L}
\newcommand{\ER}{~\!\!^E\!R}
\newcommand{\ES}{~\!\!^E\!S}

\newcommand{\OmE}{~\!\!^E\!\Om}
\newcommand{\OmS}{~\!\!^E\!\Om(\SE)}
\newcommand{\OmT}{~\!\!^E\!\Om(\Tp)}
\newcommand{\OmD}{~\!\!^E\!\Om(\Dp)}
\newcommand{\OmA}{~\!\!^E\!\Om(\Ef)}
\newcommand{\OmJ}{~\!\!^E\!\Om(\Jp)}

\newcommand{\OmEF}{~\!\!^F\!\Om}
\newcommand{\OmSF}{~\!\!^F\!\Om(\SE)}
\newcommand{\OmTF}{~\!\!^F\!\Om(\Tp)}
\newcommand{\OmDF}{~\!\!^F\!\Om(\Dp)}
\newcommand{\OmAF}{~\!\!^F\!\Om(\Ef)}
\newcommand{\OmJF}{~\!\!^F\!\Om(\Jp)}

\newcommand{\db}{\mathfrak b}
\newcommand{\pa}{\partial}
\newcommand{\mg}{\mathfrak g}

\author{Damien Calaque, Vasiliy Dolgushev and Gilles Halbout}

\title[Formality theorems for Lie
algebroid chains]{Formality theorems for Hochschild chains in
the Lie algebroid setting}

\begin{document}

\maketitle

\begin{abstract}
In this paper we prove Lie algebroid versions of
Tsygan's formality conjecture for
Hochschild chains both in the
smooth and holomorphic settings.
Our result in the holomorphic setting
implies a version of Tsygan's formality conjecture for
Hochschild chains of the structure sheaf
of any complex ma\-ni\-fold.
The proofs are based on the use
of Kontsevich's quasi-isomorphism
for Hochschild cochains of $\R[[y^1,\dots,y^d]]$,
Shoikhet's quasi-isomorphism for Hochschild chains
of $\R[[y^1,\dots,y^d]]$, and  Fedosov's resolutions of the
natural analogues of Hochschild (co)chain complexes associated with
a Lie algebroid. In the smooth setting
we discuss an application of our result
to the description of quantum traces
for a Poisson Lie algebroid.
\end{abstract}

\tableofcontents

\section*{Introduction}
Lie algebroids and Lie groupoids provide a natural framework for
developing analysis on differentiable foliations and
ma\-ni\-folds with corners  \cite{Rich}, \cite{Psi},
\cite{NWX}, \cite{W}. This motivates our interest to
the natural analogues of Hoch\-schild and cyclic (co)ho\-mo\-lo\-gi\-cal
complexes in the setting of Lie algebroids and to
the corresponding analogues of the Kontsevich-Tsygan
formality conjectures.
Thus the formality theorem for the differential graded
Lie algebra (DGLA) of Hochschild cochains in the Lie algebroid
setting \cite{C} allows us to quantize an arbitrary
Poisson Lie algebroid\footnote{According to the
terminology of P. Xu \cite{X} we have to call this object
a triangular Lie bialgebroid. However, since we do not mention 
the bialgebroid structure, we refer to this object
as a Poisson Lie algebroid.}. The formality of the DGLA module
of Hochschild chains in the Lie algebroid setting
would give a description of the quantum traces for Poisson
Lie algebroids, and the formality of the cyclic complex in
the setting of Lie algebroids would imply the algebraic index
theorem \cite{NT}, \cite{TT1} for the deformations associated with
an arbitrary Poisson Lie algebroid.

An appropriate analogue of the Hochschild cochain (resp. chain)
complex associated with a Lie algebroid $E$ is the complex of
$E$-polydifferential operators (resp. Hochschild $E$-chains)
(see definitions \ref{kotsepi} and \ref{tsepi} in the next section).
It turns out that the complex of $E$-polydifferential
operators is naturally a DGLA and the complex of
$E$-chains is naturally a DG module over this DGLA.
Due to the recent result \cite{C} of the first author
for any Lie algebroid $E$ over a smooth manifold
the DGLA of $E$-polydifferential operators is formal.

In this paper we use Kontsevich's \cite{K} and Shoikhet's \cite{S}
formality theorems for $\R^d_{formal}$ and the `Fedosov-like'
\cite{F1} globalization technique \cite{CFT,D1,D2,NT} to prove that
for any Lie algebroid $E$ over a smooth manifold
(resp. holomorphic Lie algebroid over a complex
manifold) the DGLA module of
$E$-chains (resp. the sheaf of DGLA modules of
$E$-chains) is formal. In the smooth setting
this result allows us to describe quantum
traces for an arbitrary Poisson Lie algebroid.
In the holomorphic setting this result implies
a version of Tsygan's formality conjecture
for Hochschild chains of the structure sheaf
of any complex manifold.

Eliminating the sheaf of Hochschild $E$-chains
in the holomorphic setting we get that for
any holomorphic Lie algebroid $E$ the
sheaf of $E$-polydifferential operators
is formal as a sheaf of DGLAs.
In particular, this result implies
Kontsevich's formality theorem for
complex manifolds, the proof of which
was formulated only for algebraic
varieties \cite{Ye}.

The paper is organized as follows.
In the first section we recall some basic facts about Lie algebroids
and define algebraic
structures on the complexes of $E$-polydifferential operators and
$E$-polyjets of an algebroid $E$.
We recall Kontsevich's
\cite{K} and Shoikhet's \cite{S} formality theorems for
$\R^d_{formal}$ and formulate our first
result, the formality of the module
of $E$-chains (see theorem \ref{thm:formality} on
page \pageref{thm:formality}).
The second section is devoted to the construction of the Fedosov resolutions of
the sheaves of $E$-polydifferential operators,
$E$-chains, $E$-polyvector fields and $E$-forms.
It is the most technical part of the paper.
Using these resolutions in section $3$,
we prove theorem \ref{thm:formality}.
In the same section
we apply this theorem to the description of quantum
traces of Poisson Lie algebroids.
In section $4$ we prove Tsygan's formality
conjecture for Lie algebroid chains in
the holomorphic setting (see theorem \ref{ya-Gilles}
on page \pageref{ya-Gilles}),
which, in particular, gives us the
formality theorem for Hochschild chains
of the structure sheaf of an arbitrary
complex manifold (see theorem \ref{Gilles1})
In the concluding section we mention an equivariant version of theorem
\ref{thm:formality} and raise some other questions. 

\subsection*{Notations}

We assume Einstein's convention for the summation over repeated
indices and omit the symbol $\wedge$ referring to a local basis
of exterior forms. The arrow $\brarrow$ denotes an $\Linf$-morphism
of $\Linf$-algebras, the arrow $\bbrarrow$ denotes a morphism
of $\Linf$-modules, and the notation
$$
\begin{array}{c}
\cL\\
\downarrow_{\,mod}\\
\cM
\end{array}
$$
means that $\cM$ is an $\Linf$-module
over the $\Linf$-algebra $\cL$\,. 
The abbreviation ``DGLA''
stands for ``differential graded Lie algebra'' and
the abbreviation ``DGA'' stands for
``differential graded associative algebra''. 
Throughout the paper (except section \ref{sect:holo}) 
we work over the field $\R$ of real numbers: unless otherwise specified, 
$M$ denotes a smooth real manifold, $\mO_M$ denotes the 
sheaf of real valued $C^\infty$-functions on $M$ and vector bundles are real vector bundles. 
Finally, we denote by the same symbol a vector bundle and its sheaf of sections. 

\subsection*{Acknowledgements}
We would like to thank G. Felder and A. Cattaneo for
their interest to this work.
The second author is partially supported by the NSF grant
DMS-9988796, the Grant for Support of Scientific
Schools NSh-1999.2003.2 and the grant
CRDF RM1-2545-MO-03.

\section{Algebraic structures associated with a Lie algebroid}

\subsection{Lie algebroids and associated sheaves}

Let us recall the following
\begin{dfn}\label{g-complex-alge}
A \emph{Lie algebroid} over a smooth manifold $M$
is a smooth vector bundle $E$ of finite rank
whose sheaf of sections is a sheaf of Lie algebras
equipped with a $\mO_M$-linear morphism of sheaves
of Lie algebras
$$
\rho:E\to TM\,.
$$
The $\mO_M$-module structure and the Lie algebra
structure on the sheaf $E$ are com\-pa\-tible in the following sense:
for any open subset $U\subset M$, any
function $f\in\mathcal O_M(U)$ and any
sections $u,v\in\G(U,E)$
\begin{equation}\label{eq:la}
[u,fv]=f[u,v]+\rho(u)(f)v\,.
\end{equation}
The map $\rho$ is called the \emph{anchor}.
\end{dfn}

\noindent{\bf Examples.}
1. The tangent bundle $TM$ on $M$ is
the simplest example of a Lie algebroid. The bracket is the
usual Lie bracket of vector fields and the anchor is
the identity map $\id:TM\to TM$. \\
\indent 2. More generally any involutive distribution (i.e.~regular foliation) 
$E\subset TM$ is a Lie algebroid over $M$. \\
\indent 3. A Lie algebroid over a point is simply a 
finite dimensional Lie algebra. 

\subsubsection{The sheaf of $E$-polyvector fields}

\begin{dfn}
The bundle $\ETp*$ of \emph{$E$-polyvector fields} is the
exterior algebra of the bundle $E$ with the shifted grading
\begin{equation}\label{ETp}
\ETp* = \bigoplus_{k\ge -1} \ETp{k}\,, \quad \qquad
\ETp{k} := \wedge^{k+1}E\,.
\end{equation}
\end{dfn}
It turns out that the Lie bracket $[\,,\,]$ on sections of 
$\ETp{0}=E$ can be naturally extended to a Lie bracket on sections of the whole
vector bundle $\ETp*$ of $E$-polyvectors (it was noticed in \cite{CW}).
Indeed, first, we defined a Lie bracket $[\,,\,]_{SN}$ 
on homogeneous sections of low degree as follows: 
\begin{equation}\label{S-N}
[f,g]_{SN}:=0\,,\quad  [u,f]_{SN}:=\rho(u)f\,,\quad
\textrm{and}\quad[u,v]_{SN}:=[u,v]\,.
\end{equation}
$$
\forall~ f,g\in\mO_M(U),~ u,v\in\G(U,E)
$$
Next, we extend $[\,,\,]_{SN}$ to sections of $\ETp*$ 
(i.e.~$E$-polyvector fields) by requiring the graded Leibniz rule with respect
to the $\wedge$-product: 
\begin{equation}\label{eq:SN}
[u,v\wedge w]_{SN}=[u,v]_{SN}\wedge
w+(-1)^{k(l+1)}v\wedge[u,w]_{SN}\,,
\end{equation}
$$
\forall~ u\in \G(U,\ETp{k}),~ v\in \G(U,\ETp{l}),~ w \in \G(U,\ETp*).
$$
In the simplest example $E=TM$ the Lie bracket
$[\,,\,]_{SN}$ coincides with the well known Schouten-Nijenhuis
bracket of ordinary polyvector fields. \\

\subsubsection{The sheaf of $E$-differential forms}~\\

\noindent The exterior algebra $\wedge^* E^{\vee}$ of the dual
bundle $E^{\vee}$ to $E$ is a natural candidate for the
bundle $^E\Omega^*_M$ of \emph{$E$-differential forms}. 
Sections of $^E\Omega^*_M$ (\emph{$E$-forms} for short) 
are endowed with the following 
\emph{$E$-de Rham differential}
\begin{eqnarray}
\Edif\omega(\sigma_0,\dots,\sigma_k) & := & 
\sum_i(-1)^{i}\rho(\sigma_i)\omega(\sigma_0,\dots,
\hat\sigma_i,\dots,\sigma_k) \label{E-dif} \\
& & +\sum_{i<j}(-1)^{i+j}\omega([\sigma_i,\sigma_j],
\sigma_0,\dots,\hat\sigma_i,\dots,\hat\sigma_j,\dots,\sigma_k)\,, \nonumber
\end{eqnarray}
$$
\si_i \in \G(U,E)\,.
$$

Another operation defined on $E$-forms is the contraction with
$E$-polyvector fields. For a $E$-polyvector field $u\in\G(U,\ETp k)$
we denote by $\io_u$ the contraction with $u$\,.
Using this contraction, the $E$-de Rham differential
(\ref{E-dif}), and the Cartan-Weil formula
\begin{equation}\label{C-Weil}
\Elie_u:=\Edif\circ \io_u+(-1)^k \io_u\circ\Edif
\end{equation}
we define the \emph{$E$-Lie derivative} of $E$-forms (over an open subset $U$) 
by the $E$-polyvector field $u\in \G(U,\ETp{k})$.

For our purposes it is more convenient to use the reversed grading
in the bundle of $E$-forms. Thus we denote by
\begin{equation}\label{EA}
\EA*=~\!\!^E\Omega^{-*}_M
\end{equation}
the corresponding bundle with reversed grading and
observe that
$\EA*$ is equipped with a structure of a graded module
over the sheaf of graded Lie algebras $\ETp*$
via the $E$-Lie derivative (\ref{C-Weil}).
Namely,
\begin{lem}
For any $u\in \G(U,\ETp{k})$ and $v\in \G(U,\ETp{l})$ one has
\begin{equation}\label{Elie}
\Elie_u\circ\Elie_v-(-1)^{kl}\Elie_v\circ\Elie_u=\Elie_{[u,v]_{SN}}\,.
\end{equation}
\end{lem}
\begin{proof}
First, it is immediate from the definition (\ref{C-Weil})
that for any $u\in \G(U,\ETp{k})$
\begin{equation}
\label{lie+d}
\Edif \circ \Elie_u = (-1)^k \Elie_u\circ\Edif\,.
\end{equation}
Second, we claim that
for any $v\in \G(U,\ETp{l})$
we have
\begin{equation}
\label{lie+i}
\Elie_u\circ \io_v -(-1)^{k(l+1)}\io_v\circ \Elie_u= (-1)^k
\io_{[u,v]_{SN}}\,.
\end{equation}
Using (\ref{lie+d}) and (\ref{lie+i}) it is not hard
to show that 
$$
\Elie_u (\Edif \io_v + (-)^l \io_v \Edif) -
(-)^{kl} (\Edif \io_v + (-)^l \io_v \Edif) \Elie_u =
(\Edif \io_{[u,v]_{SN}} + (-)^{k+l} \io_{[u,v]_{SN}} \Edif)\,.
$$
Thus it suffices to prove that equation (\ref{lie+i})
holds.

The proof of (\ref{lie+i}) goes as follows. First, direct
computations show that (\ref{lie+i}) holds
for any sections $u$ and $v$ of the subsheaf
$\ETp{-1}\oplus \ETp{0}$\,.
Second, using the Leibniz rule (\ref{eq:SN})
we prove the desired identity by induction on
the degrees of $E$-polyvector fields $u$ and $v$.
In doing this, we need another simple identity
$$
\Elie_{u_1\wedge u_2} = \Elie_{u_1} \io_{u_2}
- (-1)^{k_1} \io_{u_1} \Elie_{u_2}\,,
\qquad \forall~ u_i\in \G(U,\ETp{k_i})\,,
$$
which follows easily from the
fact that $\io_{u_1 \wedge u_2} = \io_{u_1} \circ \io_{u_2}$. 
\end{proof}~

\subsubsection{The sheaf of $E$-differential operators}~\\

\noindent One can also define the (left) $\mathcal O_M$-module $\Op$ of
\emph{$E$-differential operators} to be the sheaf of algebras locally generated
by functions and $E$-vector fields. More precisely,
$\Op$ is the sheaf associated with the
following presheaf
\begin{equation}\label{pre-Op}
U\longmapsto
\begin{array}{c}
T\Big(\mathcal O_M(U)\oplus\Gamma(U,E)\Big)/_{
\left\{\begin{array}{c}
f\otimes g-fg,~f\otimes u-fu, \\
u\otimes f-f\otimes u- \rho(u)f, \\
u\otimes v-v\otimes u-[u,v],
\end{array}\right\}}
\end{array}
\end{equation}
$$
f,g\in \mO_M(U), \qquad u,v \in \G(U,E)\,.
$$
As a sheaf of $\mO_M$-modules, $\Op$
is endowed with an increasing filtration
\begin{equation}\label{Op-filt}
\mO_M=\Op^0\subset \Op^1 \subset \Op^2 \subset \dots
\subset \Op\,,
\end{equation}
which is defined by assigning the degree $1$
to the $E$-polyvector fields.

In the terminology of \cite{R} $E$ is a sheaf
of \emph{Lie-Rinehart algebras} over the structure sheaf
$\mO_M$ and $\Op$ is its universal enveloping algebra. 
Besides the fact that $\Op$ is a sheaf of algebras, 
$\Op$ is also equipped with a coassociative $\mO_M$-linear map
$\Delta:\Op\to\Op\otimes_{\mathcal O_M}\Op$
which is defined as follows
$$
\Delta(f)=f\otimes1=1\otimes f\,,
$$
\begin{equation}\label{coprod}
\Delta(u)=u\otimes 1+ 1 \otimes u,
\qquad
\Delta(PQ)=\Delta(P)\Delta(Q)\,,
\end{equation}
$$
\forall~ u\in\G(U,E),~ P,Q\in \G(U,\Op)\,.
$$
{\bf Remark.}
For any $u\in\G(U,\mO_M\oplus E)$ one can see that a lift of $\Delta(u)$ lies in the normalizor $N(\mathcal I_U)$ 
of the right ideal $\mathcal I_U$ generated by $f\otimes1-1\otimes f$, $f\in\mO_M(U)$, in 
$\G(U,\Op\otimes_{\R}\Op)$. Therefore $\Delta$ takes values in a sheaf of algebras (the one associated 
to the presheaf of algebras $U\mapsto N(\mathcal I_U)/\mathcal I_U$); hence $\Delta(PQ)=\Delta(P)\Delta(Q)$ 
is well-defined. 

Moreover the anchor $\rho$ extends to a (left) $\mO_M$-linear morphism of sheaves of (associative) algebras 
$\rho:\Op\to{\rm End}(\mO_M)$. 
In the terminology of \cite{X} $(\Op,\Delta,\rho)$ is a sheaf of \emph{Hopf algebroids with anchor}. 

Notice that, in the simplest example $E=TM$ of the Lie algebroid
$\Op$ is the sheaf of differential operators on $M$\,. In this case 
$\Delta(P)$ is the bidifferential operator $(f,g)\mapsto P(fg)$. \\

The following result shows that $\Op$ is an ind-finite dimensional vector bundle over $M$. 
\begin{prop}[\cite{NWX,R}]
$\Op\cong S(E)$ as sheaves of (left) $\mO_M$-modules. 
\end{prop}

\subsection{Lie algebroids connections}

By the word \emph{connection}\index{connection} on a vector bundle $\Bu$ over $M$ we always mean $E$-\emph{connection}, 
\index{$E$-connection}that is a linear operator
\begin{equation}\label{E-conn}
\n : \G(M,\,\Bu) \to ~^E\Om^1(M,\,\Bu)
\end{equation}
satisfying the following equation
\begin{equation}\label{connexion}
\n (f u) = \Edif(f) u + f \n(u)
\end{equation}
for any $f\in \mO_M(M)$ and $u\in \G(M,\,\Bu)$.

Locally, $\nabla$ is completely determined by its \index{Christophel's symbols}\emph{Christophel's symbols} $\G_{ij}^k$. 
Namely, let $(e_1,\dots,\,e_r)$ and $(\xi^1,\dots,\,\xi^r)$ be dual local basis of $E$ and $E^\vee$ respectively, and 
$(b_1,\dots,\,b_s)$ be a local base of $\Bu$, then 
\begin{equation}\label{Chis}
\nabla(b_j)=\xi^i\G_{ij}^kb_k
\end{equation}

For any $u\in\G(M,\,E)$ we denote by $\nabla_u$ the associated map $\G(M,\,\Bu)\to\G(M,\,\Bu)$. \\
{\bf Remark.}
As with usual connections, one can extend this \emph{covariant derivative} on $E$-tensor in a unique way such that 
$\nabla_u$ is a derivation with respect to the tensor product of $E$-tensors, commutes with the contraction of 
$E$-tensors, acts as $\rho(u)$ on functions, and is $\R$-linear. 

\begin{dfn}
The \emph{curvature} \index{curvature}$R$ of a connection $\nabla$ with value in $\Bu$ is the section 
$R$ of the bundle $E^\vee\otimes E^\vee\otimes\Bu^\vee\otimes\Bu$ defined by
\begin{equation}
R(u,v)w=\big(\nabla_u\nabla_v-\nabla_v\nabla_u-\nabla_{[u,v]}\big)w
\end{equation}
for any $u,v\in\G(M,\,E)$ and $w\in\G(M,\,\Bu)$. 
\end{dfn}
Locally, the curvature is given by 
$$R(e_i,e_j)b_k=(R_{ij})_k^lb_l$$
with 
\begin{equation}\label{eq-R}
(R_{ij})_k^l=\Gamma_{im}^l\Gamma_{jk}^m-\Gamma_{ik}^m\Gamma_{jm}^l+\rho(e_i)\cdot\Gamma_{jk}^l
-\rho(e_j)\cdot\Gamma_{ik}^l-c_{ij}^m\Gamma_{mk}^l
\end{equation}

For a connection $\nabla$ on $E$ itself one has the following 
\begin{dfn}
The \emph{torsion} \index{torsion}$T$ of $\nabla$ is a $E$-tensor of type $(1,2)$ defined by 
\begin{equation}
T(u,v)=\nabla_uv-\nabla_vu-[u,v]
\end{equation}
for any $u,v\in\G(M,\,E)$. 
\end{dfn}
\noindent One can write the local coefficients of this tensor very easily: 
\begin{equation}\label{eq:T}
T_{ij}^k=\Gamma_{ij}^k-\Gamma_{ji}^k-c_{ij}^k
\end{equation}
\begin{prop}\label{torsionfree}
A torsion free connection on $E$ exists. 
\end{prop}
\begin{proof}[Proof.]
Let $(U_\alpha)_\alpha$ be a cover of $M$ by trivializing opens for $E$. On each $U_\alpha$ one has a 
base $(e_i)_i$ of sections and then can define $\nabla^{(\alpha)}_{e_i}e_j=\frac12[e_i,e_j]$. Let $(f_\alpha)_\alpha$ be 
partition of unity for $(U_\alpha)_\alpha$ and define $\nabla=f_\alpha\nabla^{(\alpha)}$. 
$\nabla$ is a torsion free connection on $E$. 
\end{proof}
\begin{prop}[Bianchi's identities]
Let $\nabla$ be connection on $E$. For any $u,\,v,\,w\in\G(M,\,E)$ one has 
\begin{equation}\label{eq:B1}
\nabla_uR(v,w)+R(T(u,v),w)+c.p.(u,v,w)=0
\end{equation}
and
\begin{equation}\label{eq:B2}
R(u,v)w-T(T(u,v),w)-\nabla_uT(v,w)+c.p.(u,v,w)=0
\end{equation}
\end{prop}
\begin{proof}[Proof.]
See for example \cite{Fer}. 
\end{proof}

\subsection{Algebraic structures on $E$-polydifferential
operators and $E$-polyjets}

\subsubsection{The sheaf of $E$-polydifferential operators}

\begin{dfn}\label{kotsepi}
The (ind-finite dimensional) graded bundle $\EDp*$ of \emph{$E$-polydifferential
operators} is the tensor algebra of the bundle
$\Op$ with a shifted grading: 
$$
\EDp* = \bigoplus_{k\ge -1} \EDp{k}\,, \quad \qquad
\EDp{k} = \otimes_{\mathcal O_M}^{k+1}\Op\,.
$$
\end{dfn}
It is easy to see that in the case $E=TM$ the
sheaf $\EDp*$ is the sheaf of polydifferential
operators on $M$.

Using the coproduct (\ref{coprod}) in $\Op$
we endow the graded sheaf
$\EDp*$ of $E$-polydifferential operators with
a Lie bracket $[,]_G$. To introduce this bracket
we first define the following bilinear
product of degree $0$
$$
\bul\,:\,\EDp{} \otimes \EDp{} \to \EDp{}\,,
$$
\begin{equation}\label{bullet}
\begin{array}{c}
\displaystyle
P\bullet Q=\sum_{i=0}^{\vert P\vert}(-1)^{i\vert Q\vert}
\big(1^{\otimes i}\otimes\Delta^{(\vert Q\vert)}
\otimes 1^{\otimes\vert P\vert-i}\big)(P)
\cdot(1^{\otimes i}\otimes Q\otimes1^{\otimes\vert
P\vert-i})\,,\\[0.4cm]
\displaystyle
P\bul f = \sum_{i=0}^{|P|}(-1)^{i}
(1^{\otimes i}\otimes \rho \otimes 1^{\otimes |P|-i} )(P) ( 1^{\otimes i}
\otimes f\otimes 1^{\otimes |P|-i})\,,\\[0.4cm]
\displaystyle
f \bul g=0\,, \qquad
f\bul P=0\,,
\end{array}
\end{equation}
for any $P, Q \in \G(U,\EDp{\ge 0})$ and $f,g\in \G(U,\EDp{-1})=\mO_M(U)$\,.
Here $\Delta^{(n)}=(\Delta\otimes 1^{\otimes
n-1})\circ\cdots\circ\Delta$, $\Delta^{(0)}$ is
by convention the identity map. \\
{\bf Remark.}
Let $\mathcal I$ be the right ideal in $\Op^{\otimes k}$ generated by 
$1^{\otimes i-1}\otimes f\otimes1^{\otimes k-i}-1^{\otimes i}\otimes f\otimes1^{\otimes k-i-1}$, 
where $i=1,\dots,k-1$ and $f\in\mO_M$. 
The sheaf of algebras $N(\mathcal I)/\mathcal I$ acts on any tensor product 
$V_{1}\otimes_{\mO_M}\cdots\otimes_{\mO_M}V_{k}$ over $\mO_M$ of left $\Op$-modules 
$V_{i}$'s. 
Since $\Delta^{(r-1)}$ obviously takes values in $N(\mathcal I)/\mathcal I$ then equation (\ref{bullet}) 
is well-defined. 

Although the bilinear product is not associative,
the graded commutator
\begin{equation}\label{Gerst}
[P,Q]_G=P\bullet Q-(-1)^{\vert P\vert\vert Q\vert}Q\bullet P\,,
\qquad P,Q\in \G(U,\EDp*)\,.
\end{equation}
defines a graded Lie bracket between the
$E$-polydifferential operators.

It is not hard to see that in the case $E=TM$
the above bracket reduces to the well known Gerstenhaber
bracket \cite{G} between polydifferential operators on $M$.

Notice that an element $1\otimes 1\in \G(M,\EDp{1})$ is
distinguished by the following remarkable
identity  $[1\otimes 1, 1\otimes 1]_G=0$\,.
Using this observation we define the following
differential
\begin{equation}\label{pa}
\pa = [1\otimes 1, ~ ]_{G}\,:\,
\EDp* \to \EDp{*+1}
\end{equation}
on the sheaf of $E$-polydifferential operators.

By definition we see that $\pa$ is compatible
with the Lie bracket (\ref{Gerst}). 
Thus, $(\EDp*,\pa,[,]_G)$ is a sheaf of
differential graded Lie algebras (DGLA for short). \\

We would like to mention that
the tensor product of sections (over $\mO_M$) turns
the sheaf $~\!\!^E\!D_{poly}[-1]^*$ with the shifted grading
into a sheaf of graded associative algebras.
Moreover, it is not hard to see that the differential $\pa$
(\ref{pa}) is compatible with this product.
Thus $\EDp*$ can be also viewed as a sheaf
of DG associative algebras (DGA). \\

\noindent{\bf Remark.}
Notice that this construction works not only for $\Op$ but for any (sheaf of) Hopf 
algebroid with anchor. 
Below, we use the fact that any morphism of Hopf algebroids with anchor induces 
a morphism between the corresponding DGLAs (resp.~DGAs). 

\subsubsection{The sheaf of $E$-polyjets}

\begin{dfn}
The bundle $\EJp*$ of \emph{$E$-polyjets}
is the following graded bundle placed in nonnegative degrees
$$
\EJp* = \bigoplus_{k\ge 0} \EJp{k}\,, \quad \qquad
\EJp{k} :=
Hom_{\mO_M}(\Op^{\otimes_{\mathcal O_M}^{k+1}}, \mO_M)\,.
$$
\end{dfn}
Since the sheaf $\EDp*$ of $E$-polydifferential operators
is an ind-finite dimensional graded vector bundle the sheaf
$\EJp*$ of $E$-polyjets is a profinite dimensional graded
vector bundle. Furthermore, the sheaf $\EJp*$ is endowed
with a canonical flat connection $\nabla^G$ which is
called the \emph{Grothendieck connection} and defined by
the formula
\begin{equation}
\label{eq:gro}
\nabla^G_\sigma(j)(P):=\rho(\sigma)(j(P))-j(\sigma \bul P)\,,
\end{equation}
where $\sigma\in\Gamma(U,E)$, $j\in\G(U,\EJp{k})$,
$P\in \G(U,\EDp{k})$, and the operation $\bul$ is
defined in (\ref{bullet}).

For this connection we have
the following standard
\begin{prop}\label{pr:chi}
Let $\chi$ be a map of sheaves
$$
\chi:\EJp{k}\to
\begin{cases}
\begin{array}{cc}
\EJp{k-1}\,, & {\rm if} \quad k > 0\,, \cr
\mO_M\,, & {\rm if} \quad k=0
\end{array}
\end{cases}
$$
defined by the formula
\begin{equation}\label{chi}
\chi(a)(P)=a(1\otimes P)\,, \quad
P \in \G(U,\EDp{k-1})\,, \quad a\in \G(U,\EJp{k})\,.
\end{equation}
The restriction of the
map $\chi$
to the $\n^G$-flat $E$-polyjets
gives the isomorphism of
sheaves
\begin{equation}\label{chi1}
\chi: \ker \n^G \cap \EJp{k}
\erarrow
\begin{cases}
\begin{array}{cc}
\EJp{k-1}\,, & {\rm if} \quad k > 0\,, \cr
\mO_M\,, & {\rm if} \quad k=0\,.
\end{array}
\end{cases}
\end{equation}
\end{prop}
\begin{proof}
To see that the map (\ref{chi1})
is surjective one has to notice that
for any $E$-polyjet $b$ of degree $k-1$
(resp. any function $b$)
the equations
$$
a(1\otimes P) = b(P)\,,  \quad P\in \G(U,\EDp{k-1})
$$
and
\begin{equation}\label{reshaiu}
a (u\cdot Q\otimes P) = \rho(u) a(Q\otimes P)
- a (Q\otimes (\D^{(k-1)}(u)\cdot P))\,,
\end{equation}
$$
Q\in \G(U,\Op)\,, \qquad u \in \G(U,E)
$$
define a $\n^G$-flat $E$-polyjet $a$ of
degree $k$ (resp. a $\n^G$-flat $E$-jet $a$)\,.

On the other hand, if $a$ is a $\n^G$-flat $E$-polyjet of
degree $k$ equation (\ref{reshaiu}) is
automatically satisfied. Thus $a$ is uniquely
determined by its image $\chi(a)$.
\end{proof}

Let $t$ be the cyclic permutation acting on
the sheaf $\EJp*$ of $E$-polyjets
\begin{equation}\label{eq:t}
t(a)(P_0\otimes\cdots\otimes P_l):=
a (P_1\otimes\cdots\otimes P_l\otimes P_0)\,,
\end{equation}
$$
a\in \G(U,\EJp{l})\,, \qquad
P_i \in \G(U,\Op)\,.
$$
Using this operation and the bilinear
product (\ref{bullet}) we define the map
\begin{eqnarray*}
\ES : \EDp{k} \otimes \EJp{l} & \to & \EJp{l-k}\,, \\
P \otimes a & \mapsto & \ES_P(a) 
\end{eqnarray*}
such that for $P\in\G(U,\EDp{k})$, $a\in\G(U,\EJp{l})$, and $Q\in \G(U,\EDp{l-k})$, 
\begin{equation}\label{ono}
\ES_P(a)(Q)=a(Q\bullet P)
+ \sum_{j=1}^k (-1)^{lj}t^j(a)
\big((\Delta^{(k)}\otimes 1^{\otimes (l-k)})(Q)\cdot
(P \otimes 1^{\otimes (l-k)})\big)\,.
\end{equation}
Due to the following proposition the map
$\ES$ defines an action
of the sheaf of graded Lie algebras $\EDp*$
of $E$-polydifferential operators
on the graded sheaf $\EJp*$
of $E$-polyjets. Namely,
\begin{prop}
\label{pr:e-moe}
For any pair $P_1, P_2\in \G(U,\EDp*)$ of
$E$-polydifferential operators
and any $E$-polyjet $a\in \G(U,\EJp*)$
\begin{equation}\label{e-moe}
\ES_{P_1} \ES_{P_2} (a) - (-1)^{|P_1| |P_2|}~\ES_{P_2}
\ES_{P_1}(a) = \ES_{[P_1, P_2]_G}(a)\,.
\end{equation}
Moreover, the action (\ref{ono}) is compatible with the
Grothendieck connection (\ref{eq:gro})
\begin{equation}
\label{ES-gro}
\n^G_{u} \left(\ES_{P_1} (a)\right) = \ES_{P_1}(\n^G_u (a))\,, \quad
u\in \G(U,E)\,.
\end{equation}
\end{prop}
\begin{proof}
It is not hard to show that 
\begin{equation}\label{Ezra}
\ES_{P_1}\ES_{P_2}(a) = \ES_{P_1\bul P_2} (a) +
H(P_1, P_2)(a) + (-1)^{|P_1| |P_2|} H(P_2, P_1)(a)\,,
\end{equation}
where\footnote{Formula (\ref{Ezra}) is essentially borrowed from
paper \cite{Ezra} of E. Getzler.}
$$
H(P_1, P_2) : \EJp{*} \to \EJp{*-|P_1|-|P_2|}
$$
is a graded $\mO_M$-linear endomorphism
of the sheaf $\EJp*$ defined by the following
formula
$$
(H(P_1,P_2)(a))(Q) =
$$
$$
\sum_{i,j} (-1)^{i|P_1| + j |P_2|}
a\left[(
1^{\otimes i} \otimes \D^{|P_1|}
\otimes 1^{\otimes(j-i-|P_1|-1)} \otimes \D^{|P_2|}
\otimes 1^{\otimes (n-j-|P_2|)}
(Q))\cdot\right.
$$
$$
\left.
(1^{\otimes i} \otimes P_1
\otimes 1^{\otimes (j-i-|P_1|-1)} \otimes P_2
\otimes 1^{\otimes (n-j-|P_2|)})
\right] +
$$
$$
\sum_{k,l} (-1)^{k|P_2|+ l(n-|P_2|)}
t^l(a)\left[ (
\D^{|P_1|} \otimes 1^{k+l-|P_1|-1}
\otimes \D^{|P_2|} \otimes 1^{\otimes n-k-l -|P_2|}
\right. (Q))\cdot
$$
$$
\left.
P_1 \otimes 1^{\otimes(k+l-|P_1|-1)}
\otimes P_2 \otimes 1^{\otimes (n-k-l-|P_2|)}
\right]\,,
$$
the sums run over all $i,j,k,l$ satisfying
the conditions
$$
0\le i \le j-|P_1|-1,\qquad  j\le n- |P_2|\,,
$$
$$
1\le l \le |P_1|, \qquad |P_1|-l+1\le k \le n-|P_2|-l\,,
$$
and
$$
Q\in \G(U,\EDp{n-|P_1|-|P_2|})\,.
$$
Equation (\ref{Ezra}) obviously implies
identity (\ref{e-moe}).

Equation (\ref{ES-gro}) follows immediately from
the fact that the coproduct (\ref{coprod}) is
compatible with the multiplication of the
$E$-differential operators and the fact that
the Grothendieck connection (\ref{eq:gro}) commutes
with the cyclic permutation (\ref{eq:t}).
\end{proof}

Like in the case of $E$-polydifferential operators the element $1\otimes 1\in \G(M,\EDp{1})$ is 
(satisfying  $[1\otimes 1, 1\otimes 1]_G=0$) allows us to define the following differential 
\begin{equation}\label{+db}
\db:=\ES_{1\otimes1}\,:\, \EJp* \to \EJp{*-1}
\end{equation}
on the sheaf of $E$-polyjets.

From the definition of the differentials (\ref{pa}), (\ref{+db})
and equation (\ref{e-moe}), we see that $\db$ is compatible with the action
(\ref{ono}) in the sense of the following equation
$$
\db \left(\ES_P (a)\right) = \ES_{\pa P} (a)
+ (-)^{|P|}\phantom{a}\ES_{P} (\db (a))\,.
$$
$$
\forall~ a\in \G(U,\EJp*)\,,\quad P\in \G(U,\EDp*)\,.
$$
Thus, $(\EJp*,\db, \ES)$ is a sheaf of differential graded modules
(DG modules for short) over $\EDp*$. \\

\subsubsection{Hochschild $E$-chains}~\\

The complex of sheaves $(\EJp*,\db)$ is not a good
candidate for the Hochschild chain complex in the Lie
algebroid setting. Indeed, if our Lie algebroid $E$
is $TM$ then the complex $(\EJp*,\db)$ boils down
to the Hochschild chain complex of $\mO_M$ without
the zeroth term and the action (\ref{ono}) does
not coincide with the standard action of
Hochschild cochains on Hochschild chains (see eq. $(3.4)$ in
\cite{D2}). To cure these problems simultaneously
we introduce a graded sheaf $\ECp*$ of
$\mO_M$-modules placed in non-positive degrees
\begin{equation}\label{ECp}
\ECp{k}= \begin{cases}
\mO_M\,, & {\rm if}~k =0\,, \\
\EJp{-k-1}\,, & {\rm if}~k< 0\,.
\end{cases}
\end{equation}
and the following $\R$-linear isomorphism
of sheaves
\begin{equation}
\label{vr}
\vr: \ECp* \to \ker \n^G\cap \EJp{-*}
\end{equation}
obtained by inverting the map (\ref{chi1}).

Due to propositions \ref{pr:e-moe}
the action (\ref{ono}) and the differential $\db$ (\ref{+db})
commute with the Grothendieck
connection $\n^G$. Thus, the $\n^G$-flat $E$-polyjets
form a sheaf of DG submodule of $(\EJp*,\db, \ES)$ over the
sheaf of DGLAs $(\EDp*, \pa, [,]_G)$\,. Combining
this observation with proposition \ref{pr:chi}
we conclude that the isomorphism (\ref{vr}) allows us
to endow the sheaf (\ref{ECp}) with a structure of
a sheaf of DG modules over the sheaf of DGLAs
$\EDp*$. Namely, 
\begin{prop}
The map
\begin{equation}
\ER_{\bul} : \EDp{k} \otimes  \ECp{l} \to \ECp{k+l}
\end{equation}
given by the formula
\begin{equation}
\label{ER}
\ER_{P}(a) = \chi \ES_{P}(\vr(a)), \quad P\in \G(U,\EDp{k}),
\quad a\in \G(U,\ECp{l})
\end{equation}
and the differential
\begin{equation}\label{db}
\db(a) = \chi \ES_{1\otimes 1} (\vr(a))
: \ECp* \to \ECp{*+1}
\end{equation}
turn $\ECp*$ (\ref{ECp}) into a sheaf of
DG modules over the sheaf of DGLAs $\EDp*$\,. $\Box$
\end{prop}
~\\
{\bf Remark 1.} Since the map $\vr$ is NOT $\mO_M$-linear the
DGLA module structure (\ref{ER}), (\ref{db}) on $\ECp*$ is only
$\R$-linear unlike the DGLA module structure (\ref{ono}) (\ref{+db})
on the sheaf $\EJp*$.\\
{\bf Remark 2.} It is not hard to see that in the
case $E=TM$ the global sections of the sheaf $\ECp*$
give the jet version \cite{Ts} of the homological
Hochschild complex of the algebra $\mO_M$ of
functions on $M$.

The second remark motivates the following
definition:
\begin{dfn}\label{tsepi}
We refer to the sheaf $\ECp*$ of DG modules
over the sheaf of DGLAs $\EDp*$ of $E$-polydifferential
operators as \emph{the sheaf of the Hochschild $E$-chains}
(or just \emph{$E$-chains} for short).
\end{dfn}

\section{The formality theorem for $E$-chains}

\subsection{Hochschild-Kostant-Rosenberg}

The cohomology of the complexes $\EDp*$ and $\ECp*$ are described
by Hochschild-Kostant-Rosenberg type theorems.
The original version of this
theorem \cite{HKR} says that the module of Hochschild homology of a smooth affine
algebra is isomorphic to the module of exterior forms of
the corresponding affine variety. In \cite{Co} A. Connes proved an
analogous statement for the algebra of smooth functions on any compact real manifold,
and in \cite{Tel}, N. Teleman was able to get rid of the assumption
of compactness. The similar question about Hochschild cohomology
turns out to be tractable if we replace the Hochschild
cochains by polydifferential operators. We believe that
the cohomology of this complex of polydifferential operators
was originally computed by J. Vey \cite{V}. All these
computations correspond to the case when $E=TM$.
In our general case we have the following proposition:
\begin{prop}\label{thm:hkr}
The natural maps
\begin{eqnarray}
\mathcal V:(\ETp*,0) & \longrightarrow & (\EDp*,\pa) \nonumber \\
v_0\wedge\cdots\wedge v_k & \longmapsto &
\frac{1}{(k+1)!}\sum_{\sigma\in S_{k+1}}\epsilon(\sigma)v_{\sigma_0}\otimes
\cdots\otimes v_{\sigma_k} \label{eq:hkr1}
\end{eqnarray}
and
\begin{eqnarray}
\mathfrak C:(\ECp*,\db) & \longrightarrow & (\EA*,0) \nonumber \\
a & \longmapsto & (v\mapsto a\circ\mathcal V(v)) \label{eq:hkr2}
\end{eqnarray}
are quasi-isomorphisms of (sheaves of) complexes.
\end{prop}
~\\
{\bf Remark.} Recall that $\EA*$ is the sheaf (\ref{EA}) of
$E$-forms with reversed grading.
\begin{proof}
It is proved in \cite{C} (Theorem 1.2) that $\mathcal V$ is a quasi-isomorphism of cochain complexes. 
By ($\mO_M$-)duality we obtain a quasi-isomorphism $(\ECp*,\flat)\to(\EA*,0)$, where $\flat a:=a\circ\partial$. 
Let us show that $\flat=\mathfrak{b}$: let $a\in\ECp{k}$ and $P\in\EDp{k-1}$, then 
\begin{eqnarray*}
(\mathfrak{b}a)(P) & = & \rho(a)((1\otimes P)\bullet m_0+(-1)^{k-1}(1\otimes P\otimes1)) \\
& = & \rho(a)(1\otimes1\otimes P-1\otimes(P\bullet m_0)+(-1)^{k-1}(1\otimes P\otimes1)) \\
& = & a(1\otimes P-P\bullet m_0+(-1)^{k-1} P\otimes1)=a(\partial P)\,.
\end{eqnarray*}
The proposition is proved. 
\end{proof}

\subsection{The formality of the DGLA module of $E$-chains}\label{subsec:main}

Unfortunately, the maps (\ref{eq:hkr1}) and (\ref{eq:hkr2}) respect neither
the Lie brackets nor the actions. This defect can be cured using the notion
of Lie algebras and their modules \emph{up to homotopy} (see \cite{HS} for a
detailed discussion of the general theory and its applications, and
\cite[section 2]{D2} for a
quick review of the notions and results we need).
The main result of this paper is the following
theorem:
\begin{thm}\label{thm:formality}
For any $C^{\infty}$ Lie algebroid
$(E,M,\rho)$ one can construct a commutative
diagram of sheaves of DGLAs and DGLA modules
over $M$
\begin{equation}
\begin{array}{ccccccc}
\ETp* & \brarrow &\cL_1 &  \brarrow  & \cL_2 & \blarrow & \EDp* \\[0.3cm]
\downarrow_{\,mod}  & ~  & \downarrow_{\,mod}& ~ &\downarrow_{\,mod}  & ~
& \downarrow_{\,mod} \\[0.3cm]
\EA* & \bbrarrow &\cM_1 &  \bblarrow  & \cM_2 & \bblarrow & \ECp*\,,
\end{array}
\label{diag-thm}
\end{equation}
in which the horizontal arrows in the upper
row are $L_\infty$-quasi-isomorphisms of sheaves of
DGLAs and the horizontal arrows in the lower
row are $L_\infty$-quasi-isomorphisms of $\Linf$-modules.
The terms ($\cL_1$, $\cL_2$,
$\cM_1$, $\cM_2$) and the quasi-isomorphisms
of the diagram (\ref{diag-thm}) are functorial for
isomorphisms of pairs $(E, \pa^E)$,
where $E$ is a $C^{\infty}$ Lie algebroid and
$\pa^E$ is a torsion free $E$-connection on
$E$.
\end{thm}
The proof of this theorem occupies the next two sections.

We would like to mention that the functoriality
of the chain of quasi-isomor\-phisms (\ref{diag-thm}) between
the pair of sheaves of DGLA modules implies the
following interesting results
\begin{cor}
Let $(E,M, \rho)$ be
a $C^{\infty}$ Lie algebroid
equipped with a smooth action
of a group $G$. If one can construct
a $G$-invariant connection $\pa^E$ on $E$
then there exists a chain of $G$-equivariant
quasi-isomorphisms between the sheaves
of DGLA modules $(\ETp*, \EA*)$ and
$(\EDp*, \ECp*)$\,.  $\Box$
\end{cor}
In particular,
\begin{cor}
If $(E,M, \rho)$ is
a $C^{\infty}$ Lie algebroid
equipped with a smooth action
of a finite or compact group $G$
then the DGLA modules $(\G(M,\ETp*)^G$, $\G(M, \EA*)^G)$
and $(\G(M,\EDp*)^G$, $\G(M,\ECp*)^G)$
are quasi-isomorphic. $\Box$
\end{cor}

\noindent{\bf Example.} Let us consider the case 
when the base manifold $M$ shrinks to a point. 
Then the Lie algebroid $E$ is a finite 
dimensional real Lie algebra $\mg$ and the 
diagram of sheaves (\ref{diag-thm}) becomes a 
diagram of (DG) Lie algebras and their modules. 
These DG Lie algebras and their modules 
can be described in geometric terms using 
a real Lie group $G$ whose Lie algebra 
is $\mg$\,. Indeed, $\EDp*$ can be identified 
with the DGLA of the left invariant polydifferential 
operators on $G$\,, $\ECp*$ is the $\EDp*$-module 
of left invariant polyjets on $G$. Similarly,
the sheaves $\ETp{*}$ and $\EA{*}$ can be 
identified with the graded Lie algebra of left 
invariant polyvector fields on $G$ and 
the graded module of left invariant 
exterior forms on $G$\,, respectively.  
In this case, our result can be derived from 
corollary $4$ in \cite{thesis} (see section $5.3$ in 
\cite{thesis}).

~\\

\noindent{\bf Remark.}
It will appear clearly in the proof that all these results remain true for 
complex Lie algebroids. 
Namely, a complex Lie algebroid on a smooth real manifold $M$ 
is a complex vector bundle of finite rank $E$ whose sheaf of sections 
is a sheaf of (complex) Lie algebras with a $\mO_M^\C$-linear morphism 
of sheaves of Lie algebras $\rho:E\to T_\C M$ satisfying the same condition 
described in formula (\ref{eq:la}). 

\subsection{Formality theorems for the Hochschild
complexes of $\R[[y^1, \dots, y^d]]$.}

In order to prove theorem \ref{thm:formality} we
construct the Fedosov resolutions of the sheaves of
DGLAs $\ETp*$ and $\EDp*$ and of the sheaves of
DGLA modules $\EA*$ and $\ECp*$. These resolutions allow us
to reduce the problem to the case of the tangent Lie algebroid
$T\R^d\to \R^d$. For the latter case the desired result follows
from the combination of Kontsevich's \cite{K} and Shoikhet's
\cite{S} formality theorems.

First, we recall the required version of Kontsevich's formality
theorem. Let $M=\R^d_{formal}$ be the formal completion of $\R^d$
at the origin. In other words we set $\mathcal O_M=\R[[y^1,\dots,y^d]]$ and
$E={\rm Der}(\mathcal O_M)$.
Let us denote by $T^*_{poly}(\R^d_{formal})$ and
$D^*_{poly}(\R^d_{formal})$ the DGLA of polyvector fields and
polydifferential operators on $\R^d_{formal}$, respectively, then

\begin{thm}[Kontsevich, \cite{K}]\label{thm:kontsevich}
There exists an $L_\infty$-quasi-isomorphism $\mathcal K$ from
$T^*_{poly}(\R^d_{formal})$ to $D^*_{poly}(\R^d_{formal})$
such that
\begin{enumerate}
\item The first structure map $\mathcal K^{[1]}$ is Vey's
quasi-isomorphism (\ref{eq:hkr1}) of complexes $\mathcal V$.
\item $\mathcal K$ is $GL_d(\R)$-equivariant.
\item If $n>1$ then for any vector fields
$v_1,\dots,v_n\in T^0_{poly}(\R^d_{formal})$
$$\mathcal K^{[n]}(v_1,\dots,v_n)=0$$
\item If $n>1$ then for any vector field $v\in T^0_{poly}(\R^d_{formal})$
linear in the coordinates $y^i$ and any polyvector fields
$\chi_2,\dots,\chi_n\in T^*_{poly}(\R^d_{formal})$
$$\mathcal K^{[n]}(v,\chi_2,\dots,\chi_n)=0.$$
\end{enumerate}
\end{thm}

We denote by
$$
A^*(\R^d_{formal}) = \R[[y^1,\dots,y^d]]\otimes\, \bigwedge (\R^d)
$$
the complex of exterior forms on $\R^d_{formal}$
with the vanishing differential and by
$$
J_*^{poly}(\R^d_{formal}) = \R[[y^1,\dots,y^d]]^{\hotimes\, (*+1)}
$$
the complex of Hochschild chains of $\R[[y^1,\dots,y^d]]$\,, where
the notation $\hotimes$ stands for the tensor product completed
in the adic topology on $\R[[y^1,\dots,y^d]]$.

Using the Lie derivative (\ref{C-Weil}) of exterior forms by
a polyvector field, we can regard $A^*(\R^d_{formal})$ as a
graded module over the graded Lie algebra
$T^*_{poly}(\R^d_{formal})$. Furthermore,
the action of Hochschild cochains on Hochschild
chains (see formula $(3.4)$ in \cite{D2}) allows us
to regard $J_*^{poly}(\R^d_{formal})$ as a DG modules
over the DGLA $D^*_{poly}(\R^d_{formal})$.
Finally, using Kontsevich's quasi-isomorphism $\mathcal K$
we get an $L_\infty$-module structure on $J_*^{poly}(\R^d_{formal})$
over $T^*_{poly}(\R^d_{formal})$. For this $L_\infty$-module,
we have the following theorem:
\begin{thm}[Shoikhet, \cite{S}]\label{thm:shoikhet}
There exists a quasi-isomorphism $\mathcal S$ of $L_\infty$-modu\-les
over $T^*_{poly}(\R^d_{formal})$ from $J_*^{poly}(\R^d_{formal})$ to
$A^*(\R^d_{formal})$ such that
\begin{enumerate}
\item The first structure map $\mathcal S^{[1]}$ is the
quasi-isomorphism of Connes (\ref{eq:hkr2})\,.
\item The structure maps of $\mathcal S$ are $GL_d(\R)$-equivariant.
\item If $n>1$ then for any vector field $v\in T^0_{poly}(\R^d_{formal})$
linear in the coordinates,
any polyvector fields
$\chi_2,\dots,\chi_n\in T^*_{poly}(\R^d_{formal})$ and any chain
$j\in J_*^{poly}(\R^d_{formal})$
$$\mathcal S^{[n]}(v,\chi_2,\dots,\chi_n;j)=0$$
\end{enumerate}
\end{thm}
~\\
{\bf Remark 1.} The third assertion of the above theorem is
proved in \cite{D2} ~(see ~theorem ~$3$).\\
{\bf Remark 2.} Hopefully, one can prove the assertions of
theorem \ref{thm:shoikhet} along the lines of
Tamarkin and Tsygan \cite{Ta,TT1,TT2}.

\section{The Fedosov resolutions}

Let, as above, $E\to M$ be a $C^\infty$ Lie algebroid with bracket
$[,]$ on sections and the anchor $\rho$.
Following \cite{D2} we introduce the formally completed
symmetric algebra bundle $\SE$ of the dual bundle $E^{\vee}$ and
bundles $\Tp,~\Dp,~\Ef,~\Jp$ naturally associated to $\SE$.
They all are pro- and/or ind-finite dimensional vector bundles. 

\begin{itemize}
\item $\SE$ is the formally completed symmetric
algebra bundle of the bundle $E^{\vee}$\,.
Local sections are given by formal power series
$$
\sum_{l=0}^\infty s_{i_1\dots i_l}(x)y^{i_1}\cdots y^{i_l}
$$
where $y^i$ are coordinates on the fibers of $E$ and
$s_{i_1\dots i_l}$ are components of a symmetric covariant
$E$-tensor.

\item $\Tp^*:=\SE\otimes\wedge^{*+1}E$ is the graded bundle of formal fiberwise
polyvector fields on $E$. Local homogeneous sections of degree $k$ are of the form
\begin{equation}
\label{sect-Tp}
\sum_{l=0}^\infty v_{i_1\dots i_l}^{j_0\dots j_k}(x)y^{i_1}\cdots y^{i_l}
\frac\pa{\pa y^{j_0}}\wedge\cdots\wedge\frac\pa{\pa y^{j_k}},
\end{equation}
where $v_{i_1\dots i_l}^{j_0\dots j_k}$ are components of an $E$-tensor with symmetric
covariant part (indices $i_1,\dots,i_l$) and antisymmetric contravariant part
(indices $j_0,\dots,j_k$).
\item $\Dp^*:=\SE\otimes T^{*+1}(SE)$ is the graded bundle of formal fiberwise
polydifferential operators on $E$ with the shifted grading.
A local homogeneous section of degree $k$ looks as follow
\begin{equation}
\label{sect-Dp}
\sum_{l=0}^\infty P_{i_1\dots i_l}^{\alpha_0\dots\alpha_k}(x)y^{i_1}\cdots y^{i_l}
\frac{\pa^{\vert\alpha_0\vert}}{\pa y^{\alpha_0}}\otimes\cdots\otimes
\frac{\pa^{\vert\alpha_k\vert}}{\pa y^{\alpha_k}},
\end{equation}
where $\alpha_s$ are multi-indices, $P_{i_1\dots i_l}^{\alpha_0\dots\alpha_k}$ are
components of an $E$-tensor with the obvious symmetry of the corresponding
indices, and
$$
\frac{\pa^{\vert\alpha_s\vert}}{\pa y^{\alpha_s}} =
\frac{\pa}{\pa y^{j_1}} \dots \frac{\pa}{\pa y^{j_{|\alpha_s|}}}
$$
for $\al_s= (j_1 \dots j_{|\alpha_s|})$\,.

\item $\Ef_*:=\SE\otimes\wedge^{-*}(E^\vee)$ is the graded bundle of
formal fiberwise differential forms on $E$ with the reversed grading.
Any local homogeneous section of degree $-k$ can be written as
\begin{equation}
\label{sect-Ef}
\sum_{l=0}^\infty \omega_{i_1\dots i_l,j_1\dots j_k}(x)y^{i_1}\cdots y^{i_l}dy^{j_1}
\wedge\cdots\wedge dy^{j_k},
\end{equation}
where $\omega_{i_1\dots i_l,j_1\dots j_k}$ are components of a
covariant $E$-tensor symmetric in indices $i_1,\dots,i_l$
and antisymmetric in indices $j_1,\dots,j_k$.
\item $\Jp_*$ is the bundle of Hochschild
chains of $\SE$ over $\mO_M$.
\begin{equation}
\label{Jp}
\Jp = \bigoplus_{k\ge 0} \Jp_k,
\qquad \Jp_k:= (\hat SE^\vee)^{\hotimes_{\mO_M} (k+1)}\,,
\end{equation}
where $\hotimes$ stands for the tensor product completed
in the adic topology. Local sections of homogeneous degree
$k$ are formal power series
\begin{equation}
\label{sect-Jp}
\sum_{\alpha_0,\dots,\alpha_k}a_{\alpha_0,\dots,\alpha_k}(x)
y_0^{\alpha_0} y_1^{\alpha_1} \cdots y_{k}^{\alpha_k}
\end{equation}
in $k+1$ copies $y_0, \dots, y_k$ of coordinates
on the fibers of $E$. Here $\alpha_s$ are multi-indices,
$a_{\alpha_0,\dots,\alpha_k}$ are
components of a tensor with an obvious symmetry in the
corresponding indices, and
$$
y_m^{\al_m} = y_m^{j_1} \dots y_m^{j_{|\alpha_m|}}
$$
for $\al_m= (j_1 \dots j_{|\alpha_m|})$\,.
\end{itemize}

For our purposes, we consider $E$-differential forms with values in
the sheaves $\SE,~\Tp,~\Dp,~\Ef,~\Jp.$ Below we list these
sheaves of $E$-forms together with the algebraic structures
they carry.\footnote{For any bundle $\mathcal B$ we will denote $\OmE(\mathcal B)$ the bundle 
$\OmE\otimes\mathcal B$ of $E$-forms with values in $\mathcal B$, and $\OmE(U,\mathcal B)$ the space 
of sections over an open subset $U\subset M$ (instead of $\G(U,\OmE\otimes\mathcal B)$). }
\begin{itemize}
\item $\OmS$ is a bundle of graded commutative algebras with grading given
by the exterior degree of $E$-forms. $\OmS$ is also filtered by
the degree of monomials in fiber coordinates $y^i$.

\item $\OmT$ is a sheaf of graded Lie algebras and $\OmA$ is
a sheaf of graded modules over $\OmT$. These structures
are induced by those of $T_{poly}^*(\R_{formal}^d)$ and
$A^*(\R^d_{formal})$, respectively and the grading
is given by the sum of the exterior degree and
the degree of a polyvector (resp. a form).
$[,]_{SN}$ will denote the Lie bracket between sections
of the sheaf $\OmT$ and $L_{u}$ (the Lie derivative)
will denote the action of a fiberwise polyvector $u\in \OmT$
on the sections of $\OmA$. $\OmT$ is also a sheaf of
graded commutative algebras. The multiplication of
sections in $\OmT$ is given by the exterior product in the space
$T_{poly}^*(\R_{formal}^d)$ of polyvector fields
on $\R^d_{formal}$\,. The Lie bracket and the product
in $\OmT$ turn $\OmT$ into a sheaf of Gerstenhaber
algebras\footnote{The definition of the Gerstenhaber
algebra can be found in section $4.1$
of the second part of \cite{Ts-book} or in the
original paper \cite{G}.}.

\item $\OmD$ is a sheaf of DGLAs and $\OmJ$ is
a sheaf of DGLA modules over $\OmD$. These structures
are induced by those of $D_{poly}^*(\R_{formal}^d)$ and
$J^*_{poly}(\R^d_{formal})$, respectively and the grading
is given by the sum of the exterior degree and
the degree of a (co)chain. We denote by $\pa$ and $[,]_G$
respectively the differential and the Lie bracket on $\OmD$,
$\db$ will stand for the differential on $\OmJ$ and
$\cR_{P}$ will denote the action of $P\in \OmD$ on
the sections of $\OmJ$. $\OmD$ is also a sheaf of
DGAs. The multiplication of
sections is induced by the cup product in
the space $D_{poly}^*(\R_{formal}^d)$ of polydifferential
operators on $\R_{formal}^d$\,.

\end{itemize}
{\bf Remark.}
Notice that $\Ef$ is a sheaf of exterior forms
with values in $\SE$. However, we would like
to distinguish $\Ef$ from $\OmS$. For this
purpose we use two copies of a local basis
of exterior forms. Those are $\{ dy^i \}$ and
$\{ \xi^i \}$ for $\Ef$ and $\OmS$, respectively.

\medskip

The following proposition shows that
we have a distinguished sheaf of graded
Lie algebras which acts on the
sheaves $\OmS$, $\OmA$, $\OmT$, $\OmD$,
and $\OmJ$.
\begin{prop}
\label{ono1}
The sheaf $\OmE(\Tp^0)$ of $E$-forms with values
in fiberwise vector fields is a sheaf of graded Lie
algebras. The sheaves
$\OmS$, $\OmA$, $\OmT$, $\OmD$, and $\OmJ$
are sheaves of modules over $\OmE(\Tp^0)$ and
the action of sections in $\OmE(\Tp^0)$ is
compatible with the DG algebraic structures
on $\OmS$, $\OmA$, $\OmT$, $\OmD$, and $\OmJ$\,.
\end{prop}
\begin{proof}
Since the Schouten-Nijenhuis bracket (\ref{S-N}), (\ref{eq:SN})
has degree zero $\OmE(\Tp^0)$ $\subset$ $\OmT$
$\subset$ $\OmD$ is a
subsheaf of graded Lie algebras.
While the action of $\OmE(\Tp^0)$ on the
sections of $\OmS$ is obvious, the action on
$\OmA$ is given by the Lie derivative, the
action on $\OmT$ is the adjoint action corresponding
to the Schouten-Nijenhuis bracket, the action on
$\OmD$ is given by the Gerstenhaber bracket and
the action on $\OmJ$ is induced by the
action of Hochschild cochains on
Hochschild chains (see formula $3.4$ in
paper \cite{D2}). The compatibility of the action
with the corresponding DGLA and DGLA module structures follows
from the construction. The compatibility of the action with the
product in $\OmT$ follows from the axioms of
the Gerstenhaber algebra \cite{G}
and the compatibility with the
product in $\OmD$ can be verified by
a straightforward computation.
\end{proof}

Due to the above proposition
the following $2$-nilpotent derivation
\begin{equation}\label{de}
\delta := \xi^i\frac{\pa}{\pa y^i}:
\OmE^*(\SE) \to\OmE^{*+1}(\SE)
\end{equation}
of the sheaf of algebras $\OmS$
obviously extends to $2$-nilpotent
differentials on $\OmT, \OmD, \OmA$ and $\OmJ$.
Furthermore, it follows from proposition \ref{ono1}
that $\de$ is compatible with the DG algebraic
structures on $\OmT$, $\OmA$,
$\OmD$, and $\OmJ$.

Note that
\begin{equation}\label{de-SE}
\ker \de \cap \SE \cong \mO_M\,, \qquad
\ker \de \cap \Ef_* \cong \EA*
\end{equation}
as sheaves of (graded) commutative
algebras over $\mO_M$.
Similarly, $\ker \de \cap \Tp$, (resp.
$\ker \de \cap \Dp$) is a sheaf
of fiberwise polyvector fields (\ref{sect-Tp})
(resp. fiberwise polydifferential operators
(\ref{sect-Dp}))
whose components do not depend on the
fiber coordinates $y^i$. In other words,
\begin{equation}
\label{de-Tp}
\ker \de \cap \Tp^* \cong \wedge^{*+1}(E)
\end{equation}
as sheaves of graded commutative algebras and
\begin{equation}
\label{de-Dp}
\ker \de \cap \Dp^* \cong \otimes^{*+1}(S(E))\,,
\end{equation}
as sheaves of DGAs over $\mO_M$.

In fact, one can prove a
more stronger statement:
\begin{prop}\label{thm:delta}
For $\Bu$ being either of the sheaves
$\SE$, $\Ef$, $\Tp$ or $\Dp$
$$
H^{\ge 1}(\OmE(\Bu), \de) =0\,.
$$
Furthermore,
\begin{equation}\label{H-de0}
\begin{array}{c}
H^0(\OmS, \de) \cong \mO_M\,, \\[0.3cm]
H^0(\OmE(\Ef_*), \de) \cong \EA*\,, \\[0.3cm]
H^0(\OmE(\Tp^*), \de) \cong \wedge^{*+1}(E)
\end{array}
\end{equation}
as sheaves of (graded) commutative algebras
and
\begin{equation}
\label{H-de01} \qquad
H^0(\OmE(\Dp^*), \de) \cong \otimes^{*+1}(S(E))
\end{equation}
as sheaves of DGAs
over $\mO_M$.
\end{prop}
\begin{proof}
Due to equations (\ref{de-SE}), (\ref{de-Tp}),
and (\ref{de-Dp}) the proposition will follow immediately
if we construct an operator
\begin{equation}
\ka: \OmE^*(\Bu) \to \OmE^{*-1}(\Bu)
\label{kappa}
\end{equation}
such that for any section $u$
of $\OmE(\Bu)$
\begin{equation}\label{eq:homotopy}
u=\delta\kappa (u)+\kappa\delta (u) + \cH(u)\,,
\end{equation}
where
\begin{equation}
\label{cH}
\cH(u) = u\Big|_{y^i=\xi^i=0}\,.
\end{equation}
First, we define this operator on
the sheaf $\OmS$
\begin{equation}
\kappa(a) = y^k
\frac {\vec{\partial}} {\partial \xi^k}
\int\limits_0^1 a(x,t y,t\xi)\frac{dt} t,
\qquad a \in \OmE^{>0}(U,\SE)\,,
\qquad \ka\Big|_{\SE}=0\,,
\label{kappa1}
\end{equation}
where the arrow over $\pa$ denotes
the left derivative with respect to
the anti-com\-muting variable $\xi^i$.

Next, we extend $\ka$ to sections of the
sheaves $\OmA$, $\OmT$, $\OmD$
in the componentwise manner.
A direct computation shows that
equation (\ref{eq:homotopy}) holds and the
proposition follows.
\end{proof}

Since our Lie algebroid $E$ is a smooth
bundle over $M$, it admits a global torsion free
connection $\pa^E$\footnote{Recall that
by the word ``connection'' we always mean
an $E$-connection (\ref{E-conn}).}.
Using this connection
we define the following derivation of
the DG sheaves $\OmS$, $\OmA$, $\OmT$,
$\OmD$, and $\OmJ$:
\begin{equation}
\label{nabla}
\n = \Edif + \G\cdot :\OmE^{*}(\Bu) \to \OmE^{*+1}(\Bu)\,, \qquad
\G = -\xi^i \Gamma_{ij}^k y^j\frac{\pa}{\pa y^k}\,,
\end{equation}
where $\Bu$ is either of the sheaves $\SE$, $\Ef$,
$\Tp$, $\Dp$, or $\Jp$, $\Gamma_{ij}^k(x)$ are
Christoffel's symbols of the connection $\pa^E$ and
$\G\cdot$ denotes the action of $\G$ on the
sections of the sheaves $\OmE(\Bu)$ (see proposition \ref{ono1}).
It is not hard to see that $\n$ (\ref{nabla}) is
compatible with the DG algebraic
structures on $\OmS$, $\OmT$, $\OmA$,
$\OmD$, and $\OmJ$. Furthermore, the torsion
freeness of the connection $\pa^E$ implies that
\begin{equation}
\label{torsion}
\n \de +\de \n = 0\,.
\end{equation}

The standard curvature $E$-tensor $(R_{ij})_k^l(x)$ of
the connection $\pa^E$ provides us with
the following fiberwise vector
field:
\begin{equation}\label{eq:R}
R=- \frac12\xi^i \xi^j (R_{ij})_k^l(x) y^k\frac{\pa}{\pa y^l}
\in \OmE^2(M,\Tp^0)\,.
\end{equation}
A direct computation shows that for
$\Bu$ being any of the sheaves
$\ES$, $\Ef$, $\Tp$, $\Dp$, or $\Jp$, we have
\begin{equation}
\label{eq:R1}
\n^2 = R\cdot : \OmE^*(\Bu) \to \OmE^{*+2}(\Bu)\,,
\end{equation}
where $R\cdot$ denotes the action of the
vector field $R$ in the sense of
proposition \ref{ono1}.

Although $\nabla$ is not flat the following theorem
shows that the combination $\n- \de$ can be
extended to a flat connection on
the sheaves $\OmS$, $\OmT$, $\OmA$,
$\OmD$, and $\OmJ$.
\begin{thm}
\label{G-Fed}
Let $\Bu$ be either of the sheaves
$\ES$, $\Ef$, $\Tp$, $\Dp$, or $\Jp$\,.
There exists a global section
\begin{equation}
\label{eq:A}
A=\sum_{s=2}^\infty\xi^k
A_{k,i_1\dots i_s}^j(x) y^{i_1}\cdots y^{i_s}\frac{\pa}{\pa y^j}
\end{equation}
of the sheaf $\OmE^1(\Tp^0)$
such that the derivation
\begin{equation}
\label{DDD}
D:=\nabla-\delta+A\cdot : \OmE^{*}(\Bu)
\to \OmE^{*+1}(\Bu)
\end{equation}
is $2$-nilpotent
$$
D^2 = 0\,,
$$
and (\ref{DDD}) is compatible with the
DG  algebraic structure on $\OmE(\Bu)$\,.
\end{thm}
\begin{proof}
The proof goes essentially along the lines of
\cite[theorem 2]{D1}.

Thanks to equation (\ref{eq:R1}) the condition $D^2=0$ is
equivalent to the equation
\begin{equation}
\label{DDD-nil}
R + \n A - \de A  + \frac12 [A,A]_{SN} =0\,.
\end{equation}
We claim that a solution of (\ref{DDD-nil})
can be obtained by iterations of the
following equation
\begin{equation}
\label{iter}
A=\kappa R + \kappa(\nabla A+\frac12[A, A]_{SN})
\end{equation}
in degrees in the fiber coordinates $y^i$\,.
Indeed, equation (\ref{eq:homotopy}) implies
that iterating (\ref{iter}) we get
a solution of the equation
$$
\ka(R + \n A - \de A + \frac12[A,A]_{SN})=0\,.
$$

We denote by $C$ the
left hand side of (\ref{DDD-nil})
$$
C = R + \n A - \de A + \frac12[A,A]_{SN}\,,
$$
and mention that due to Bianchi's identities
$\nabla R=\delta R=0$
\begin{equation}\label{urav}
\n C - \de C + [A,C]=0\,.
\end{equation}
Applying $\ka$ (\ref{kappa1}) to (\ref{urav})
and using the homotopy property (\ref{eq:homotopy})
we get
$$
C = \ka (\n C + [A,C])\,.
$$
The latter equation has the unique vanishing
solution since the operator $\ka$ (\ref{kappa1})
raises the degree in the fiber coordinates $y^i$\,.

Proposition \ref{ono1} implies that the
differential (\ref{DDD}) is compatible with the
DG algebraic structures on $\OmE(\Bu)$\,.
Thus, the theorem is proved.
\end{proof}
In what follows we refer to the differential $D$
(\ref{DDD}) as \emph{the Fedosov differential}.

The following theorem describes the cohomology
of the Fedosov differential $D$ for the
sheaves $\OmS$, $\OmA$, $\OmT$, and $\OmD$
\begin{thm}\label{thm:resolution1}
For $\bB$ being either of the sheaves
$\OmS$, $\OmA$, $\OmT$, or $\OmD$
\begin{equation}\label{H-D>0}
H^{\ge 1}(\bB, D) =0\,.
\end{equation}
Furthermore,
\begin{equation}
\label{H-D0}
\begin{array}{c}
H^0(\OmS, D) \cong \mO_M\,, \\[0.3cm]
H^0(\OmE(\Ef_*), D) \cong \EA*\,, \\[0.3cm]
H^0(\OmE(\Tp^*), D ) \cong  \ker \de \cap \Tp^*\,,
\end{array}
\end{equation}
as sheaves of graded commutative algebras
\begin{equation}
\label{H-D01}
H^0(\OmE(\Dp^*), D) \cong   \ker \de \cap \Dp^*
\end{equation}
as sheaves of DGAs over $\R$.
\end{thm}
\begin{proof}
The first statement follows easily from the spectral sequence
argument. Indeed, using the fiber coordinates $y^i$ we
introduce the decreasing filtration
$$
\dots \subset F^{p+1} \bB \subset F^p \bB \subset F^{p-1}\bB
\subset \dots \subset F^0\bB=\bB\,,
$$
where the components of the sections of
the sheaf $F^p\bB$ have degree in $y^i$
$\ge p$.

Since $D (F^p\bB) \subset F^{p-1}\bB$ the corresponding
spectral sequence starts with
$$
E_{-1}^{p,q} = F^p \bB^{p+q}\,.
$$
It is easy to see that
$$
d_{-1} = \de\,.
$$
Thus using proposition \ref{thm:delta} we conclude that
for any $p,q$ satisfying the condition $p+q > 0$
$$
E_0^{p,q} = E_1^{p,q} = \dots = E_{\infty}^{p,q}=0
$$
and the first statement (\ref{H-D>0}) follows.

Let $\Bu$ denote either of the bundles $\SE$, $\Ef$,
$\Tp$, or $\Dp$. We claim that iterating the
equation
\begin{equation}\label{iter-u}
\la(u) = u + \ka (\n \la(u) + A\cdot \la(u))\,,
\qquad u\in \G(U,\Bu\cap \ker \de)
\end{equation}
we get a map of sheaves of graded vector
spaces
\begin{equation}
\label{lift}
\la : \Bu\cap \ker \de \to \Bu\cap \ker D\,.
\end{equation}
Here $A\cdot$ denotes the action of the
fiberwise vector field $A$, defined in
proposition \ref{ono1}.
Indeed, let $u$ be a section of $\Bu$. Then, due
to formula (\ref{eq:homotopy}) $\la(u)$ satisfies the
following equation
\begin{equation}\label{ka-Y}
\ka (D (\la(u))) =0\,.
\end{equation}
Let us denote $D \la(u)$ by $Y$
$$
Y= D \la(u)\,.
$$
The equation $D^2=0$ implies that
$$
D Y = 0
$$
which is equivalent to
\begin{equation}
\label{del-Y}
\de Y = \n Y + A \cdot Y
\end{equation}

Applying (\ref{eq:homotopy}) to $Y$ and
using equations (\ref{ka-Y}), (\ref{del-Y})
we get
$$
Y = \ka (\n Y + A\cdot Y)\,.
$$
The latter equation has the unique vanishing
solution since the operator $\ka$ (\ref{kappa1})
raises the degree in the fiber coordinates $y^i$.

The map (\ref{lift}) is obviously
injective. To prove that the map is
surjective we notice that $\cH$
$$
\cH\,:\, \Bu \to \Bu\cap \ker \de
$$
is a left inverse of the map (\ref{lift}).
Thus it suffices to prove that if
$a\in \G(U,\Bu\cap \ker D)$
and
\begin{equation}
\label{cH-a}
\cH a = 0
\end{equation}
then $a$ vanishes.

The condition $a\in \ker D$ is
equivalent to the equation
$$
\de a = \n a + A \cdot a\,.
$$
Hence, applying (\ref{eq:homotopy}) to
$a$ and using (\ref{cH-a})
we get
$$
a = \ka (\n a + A\cdot a)\,.
$$
The latter equation has the unique vanishing
solution since the operator $\ka$ (\ref{kappa1})
raises the degree in the fiber coordinates $y^i$.
Thus, the map (\ref{lift}) is bijective and
the map $\cH$
\begin{equation}\label{cHH1}
\cH\,:\, \Bu\cap \ker D \to \Bu\cap \ker \de
\end{equation}
is the inverse of (\ref{lift}).

It remains to prove that the map (\ref{lift})
is compatible with the multiplication of
the sections of the sheaf $\Bu$, where
$\Bu$ is either $\SE$, $\Ef$, $\Tp$, or
$\Dp$\,. The latter
follows immediately from the fact that
the inverse map $\cH$
\begin{equation}\label{cH-A}
\cH\,:\, \Bu \to \Bu \cap \ker \de
\end{equation}
respects the corresponding algebra
structures on $\SE$, $\Ef$, $\Tp$,
and the DGA structure
on $\Dp$\,.
\end{proof}

Let us now mention that since the Fedosov differential
(\ref{DDD}) is compatible with the graded algebraic
structures on the sheaves $\OmT$ and $\OmA$ we conclude
that $H^*(\OmT,D)$ is a sheaf of graded Lie algebras and
$H^*(\OmA, D)$ is a sheaf of graded modules over
$H^*(\OmT, D)$\,. On the other hand the above theorem
tells us that
$$
H^*(\OmA, D) = \EA*\,,
$$
and
$$
H^*(\OmE(\Tp), D) = \Tp^* \cap \ker \de,
$$
Furthermore, the sheaf $\Tp^* \cap \ker \de$
in the right hand side of the latter
equation can be canonically identified
with $\ETp*=\wedge^{*+1} E$ as a
sheaf of vector spaces.

Thus, it is natural to ask whether the graded
algebraic structures on the sheaves
$\Tp^* \cap \ker \de$ and $\EA*$ coincide
with the ones given by Lie bracket (\ref{S-N})
(\ref{eq:SN})
and the Lie derivative (\ref{C-Weil}).
A positive answer to this question is given by
the following proposition:
\begin{prop}\label{resol-T-A}
The composition
\begin{equation}\label{cH'}
\cH'=\nu \circ \cH :
 \Tp^* \cap \ker D \to \ETp*
\end{equation}
of the map
\begin{equation}\label{cHHH}
\cH\,:\, \Tp^* \cap \ker D \to \Tp^* \cap \ker \de
\end{equation}
with the identification
of the sheaves $\Tp^* \cap \ker \de$ and
$\ETp* \cong \wedge^{*+1} E$
\begin{equation}\label{nu}
\nu : \Tp^* \cap \ker \de
\erarrow \ETp*
\end{equation}
induces an isomorphism of the sheaves of
graded Lie algebras $H^*(\OmT, D)\cong \ETp*$.
The map
\begin{equation}\label{cH-A1}
\cH\,:\, \Ef_* \cap \ker D \to \EA*
\end{equation}
induces an isomorphism of the sheaves of
graded modules $H^*(\OmA, D) \cong \EA*$
over the sheaf of graded Lie algebras
$H^*(\OmT, D) \cong \ETp*$\,.
\end{prop}
\begin{proof}
The first part of the proposition is
proved in \cite{C} (see proposition $2.4$).
To prove the second part, we first remark that
the maps $\cH$ and $\nu$ are compatible with
the cup products.

Next, we show that for any $D$-closed
fiberwise differential form $\omega\in \G(U,\Ef)$ one
has
$$
\cH (d^f\omega) = \Edif\mathcal H(\omega)\,,
$$
where $d^f=dy^i \frac{\pa}{\pa y^i}$ is
the fiberwise De Rham differential
on $\Ef$\,.
Since
$$
\cH : \Ef_* \to \EA*
$$
is a
morphism of graded
commutative algebras, it is
sufficient to prove it for functions and $1$-forms:
\begin{itemize}
\item \emph{First case.}
Let $f$ be a function and
$$
\omega=\la(f)\,.
$$
A direct computation shows that
$$
\la(f)=f+y^i\rho(e_i)f~\textrm{ mod }~\vert y\vert^2\,.
$$
Therefore $d^f\omega = \rho(e_i)f dy^i~\textrm{ mod }~\vert y\vert$,
and hence, $\cH(d^f \omega)=\Edif f$.
\item \emph{Second case.} Let $\alpha=\alpha_i(x)dy^i$ be a $E$-$1$-form and
$$
\omega=\la(\alpha)\,.
$$
It is not hard to show that
$$
\la(\alpha)
 =\alpha+y^i\big(\rho(e_i)\alpha_j-
\Gamma_{ij}^k\alpha_k\big)dy^j~\textrm{ mod }~\vert y\vert^2\,.
$$
Therefore,
$$
d^f \omega=\big(\rho(e_i)\alpha_j-\Gamma_{ij}^k\alpha_k\big)dy^i\wedge dy^j
~\textrm{ mod }~\vert y\vert =
\big(\rho(e_i)\alpha_j-
\frac12c_{ij}^k\alpha_k\big)dy^i\wedge dy^j
~\textrm{ mod }~\vert y\vert\,,
$$
and hence,
$$
\cH (d^f \omega) =\Edif\alpha\,.
$$
\end{itemize}
To finish the proof we notice that for
any fiberwise polyvector field $u\in \G(U,\Tp^*)$ and any
fiberwise differential form $\omega\in \G(U,\Ef)$, the equation
$$
\mathcal H(\io_u\omega)=\io_{\mathcal H(u)}\circ\mathcal H(\omega)
$$
is obviously satisfied. The latter implies that
for any pair of $D$-closed sections $u\in \G(U,\Tp^*)$,
$\omega\in \G(U,\Ef_*)$
$$
\mathcal H(L_u\omega)=\Elie_{\mathcal H(u)}\circ\mathcal H(\omega)\,,
$$
and the proposition follows.
\end{proof}
\noindent{\bf Remark.} Actually, we have proved a slightly
stronger statement. Namely, we shown that
the maps (\ref{cH-A1}) and (\ref{cH'}) induce
an isomorphism of the sheaves of
\emph{calculi}.
$$
(H^*(\OmT, D), H^*(\OmA, D)) \cong
(\ETp*, \EA*)\,.
$$
The precise definition of
the calculi can be found in
section $4.3$ of the second part
of \cite{Ts-book}.

Let us now recall that $\Tp^0$ is
a sheaf of Lie-Rinehart algebras \cite{R} over the sheaf of algebras
$\Tp^{-1}=\SE$, and $\Dp^0$ is the universal enveloping
algebra\footnote{More precisely, $\Dp^0$ is the sheaf associated to the corresponding presheaf 
of universal enveloping algebras, like $\Op$ (\ref{pre-Op}) for the sheaf of Lie-Rinehart algebras 
$E$ over $\mO_M$. } of $\Tp^0$. 
Therefore, the inverse $(\cH')^{-1}$ of
the map (\ref{cH'}) induces
the morphism
\begin{equation}\label{map-mu}
\mu : \Op \to \Dp^0
\end{equation}
of the sheaves of Hopf algebroids with anchor, and for any
$P\in \G(U,\Op)$
\begin{equation}\label{map-mu-D}
D (\mu (P)) = 0\,.
\end{equation}
We claim that
\begin{prop}\label{mu=iso}
The map (\ref{map-mu})
gives the isomorphism
\begin{equation}\label{map-mu-mu}
\mu : \Op \to \Dp^0\cap \ker D\,.
\end{equation}
of the sheaves of Hopf algebroids with anchor.
\end{prop}
\begin{proof}
Notice that $\Op$ and $\Dp^0$ are both filtered sheaves
of algebras. The filtration on $\Op$ is defined in
(\ref{Op-filt}) and the filtration
on $\Dp^0$ is given by the degree of differential
operators.

Thanks to the results of \cite{NWX} and \cite{R}
we have the PBW theorem for Lie algebroids. This theorem
says that the associated graded module of the filtration
(\ref{Op-filt}) on $\Op$ is
$$
Gr (\Op) = S(E)
$$
the symmetric algebra of the
bundle $E$\,.

Furthermore, it is not hard to see that the
map $\mu$ is compatible with the filtrations
on $\Op$ and $\Dp^0$ and due to
theorem \ref{thm:resolution1} and
proposition \ref{thm:delta}
$\mu$ induces the isomorphism
$$
S(E) \cong \Dp^0 \cap \ker D
$$
of the associated graded sheaves
of vector spaces.
Therefore, the snake lemma
argument implies that the map (\ref{map-mu-mu}) is also an isomorphism
onto the sheaf $\Dp^0 \cap \ker D$ of $D$-flat sections
of $\Dp^0$.
\end{proof}

Let us recall that $\EDp*$ (resp. $\Dp^*$)
is the tensor algebras of $\Op$ over
$\mO_M$ (resp. the tensor algebra
of $\Dp^0$ over $\SE$).
Using this fact we extend
(\ref{map-mu}) to the
morphism
\begin{equation}\label{map-mu1}
\mu' : \EDp* \to \Dp^*\,.
\end{equation}
of sheaves of DGAs
(over $\R$) by setting
$$
\mu' \Big|_{\EDp{0}} = \mu\,, \qquad
\mu' \Big|_{\mO_M} = \la\,,
$$
where $\la$ is defined in (\ref{lift})\,.

Let us also observe that
since the map (\ref{map-mu}) is a morphism of
the sheaves of Hopf algebroids with anchor then the map (\ref{map-mu1})
a morphism of the sheaves of DGLAs (over $\R$).
Furthermore,
theorem \ref{thm:resolution1} implies that
the sheaf of DGAs $\Dp^*\cap \ker D$
is generated by the sheaf $\Dp^0 \cap \ker D$
over the sheaf of commutative
algebras $\SE\cap \ker D \cong \mO_M$.
Therefore using proposition \ref{mu=iso} we get the
following result:
\begin{prop}[proposition 2.5, \cite{C}]
\label{mu-EDp}
The map (\ref{map-mu1}) gives an isomorphism
of the sheaves of DGLAs
\begin{equation}
\label{map-mu11}
\mu' : \EDp* \erarrow \Dp^* \cap \ker D\,.
\end{equation}
This map is also compatible with the DGA
structures on the sheaves $\EDp*$ and
$\Dp^* \cap \ker D$ by construction. $\Box$
\end{prop}

Let us consider the map of
sheaves of graded vector spaces
\begin{equation}
\label{map-ga}
\ga : \Jp_* \to \EJp*\,,
\qquad
\ga(j)(P)= (\mu' (P))(j)\Big|_{y^i=0}\,,
\end{equation}
$$
j\in \G(U,\Jp_k)\,, \quad P\in \G(U,\EDp{k})\,.
$$
We claim that
\begin{thm}\label{thm-ona}
For any $q\ge 1$
\begin{equation}\label{H0D-J}
H^{q} (\OmJ,D) = 0\,,
\end{equation}
and the map (\ref{map-ga}) gives
an isomorphism of the sheaves
of DG modules over
the sheaf of DGLAs
$\EDp*\cong \Dp^*\cap \ker D$
\begin{equation}
\label{map-ga1}
\ga : \Jp_* \erarrow \EJp*\,.
\end{equation}
This isomorphism sends the Fedosov
connection (\ref{DDD}) on $\Jp^*$ to
the Grothen\-dieck connection
(\ref{eq:gro}) on $\EJp*$\,.
\end{thm}
\begin{proof}
The first statement (\ref{H0D-J}) follows easily from the spectral sequence
argument. Indeed, using the zeroth collection of the fiber coordinates
$y_0^i$ (\ref{sect-Jp}) we introduce the decreasing filtration
on the sheaf $\OmJ$
$$
\dots \subset F^{p+1} (\OmJ) \subset F^p (\OmJ) \subset
F^{p-1}(\OmJ)
\subset \dots \subset F^0(\OmJ)=\OmJ\,,
$$
where the components of the sections (\ref{sect-Jp}) of
the sheaf $F^p(\OmJ)$ have degree in $y_0^i$
$\ge p$.

Since $D (F^p(\OmJ)) \subset F^{p-1}(\OmJ)$ the corresponding
spectral sequence starts with
$$
E_{-1}^{p,q} = F^p (\OmJ^{p+q})\,.
$$
Next, we observe that
$$
d_{-1} = \xi^i \frac{\pa}{\pa y_0^i}\,,
$$
and hence, due to the Poincar\'e lemma for the formal disk
we have
$$
E_0^{p,q} = E_1^{p,q} = \dots = E_{\infty}^{p,q} =0
$$
whenever $p+q> 0$.
Thus, the first statement (\ref{H0D-J}) of the
theorem follows.

Since (\ref{map-mu}) is a morphism of sheaves of
Hopf algebroids with anchor
$$
\mu'(P\,\bul\, Q) = \mu'(P)\, \bul \, \mu'(Q)\,,
\qquad P,\,Q\in \G(U,\EDp*)\,.
$$
Furthermore, $\mu'$ is obviously compatible
with cyclic permutations
$$
t \,\mu'(P_0 \otimes P_1\otimes \dots \otimes P_l)
= \mu'(P_1 \otimes P_2\otimes \dots\otimes P_l \otimes P_0)\,,
\qquad
P_i\in \G(U,\Op)\,.
$$
Hence, for any $P\in \G(U,\EDp*)$
and any $a\in \G(U,\Jp_*)$
\begin{equation}
\label{map-ga-R}
\ES_P (\ga(a)) = \ga(\cR_{\mu'(P)}(a))\,.
\end{equation}

Since $\Jp_*$ is dual to $\Dp^*\cap \ker \de$
and $\Dp^*\cap \ker \de\cong \Dp^*\cap \ker D
\cong \EDp*$ the map (\ref{map-ga1}) is an
isomorphism.
It remains to prove that the map (\ref{map-ga1})
sends the Fedosov connection (\ref{DDD}) to the Grothendieck
connection (\ref{eq:gro}). This statement is proved
by the following line of equations:
\begin{align*}
\ga(D_u j)(P)& = (\mu'(P))(D_u j)\Big|_{y^i=0} =
(D_u [\mu'(P)(j)]) \Big|_{y^i=0}\cr
&=\rho(u)[\mu'(P)(j)]\Big|_{y^i=0}
- (\io_u \de \, \bul \, [\mu'(P)(j)] )\Big|_{y^i=0}\cr
&=
\rho(u) [\mu'(P)(j)]\Big|_{y^i=0}
- (\mu'(u)\, \bul \, \mu'(P))(j)\Big|_{y^i=0}
\cr
&=\rho(u) [\mu'(P)(j)]\Big|_{y^i=0}
- \mu'(u \, \bul \,  P)(j)\Big|_{y^i=0}
\cr
&=\rho(u) (\ga(j))(P) - (\ga(j))(u \, \bul \,  P) = (\n^G_u \ga(j))(P)\,,
\end{align*}
where $u\in \G(U,E)$, $j \in \G(U,\Jp_{k})$,
$P\in \G(U,\EDp{k})$\,, $\io$ denotes the
contraction of an $E$-vector field with
$E$-differential forms, $\rho$ is the
anchor map, and $u$ is viewed
both as a section of $E$ and an $E$-differential
operator.
\end{proof}

\section{Proof of the formality theorem for $E$-chains and
its applications}

\subsection{Proof of the theorem}

Let us denote
\begin{itemize}
\item $\la_{A}$ : $\EA*\to
\OmE(\Ef_*)$ the map $\la$ (\ref{lift}) defined in the proof of theorem \ref{thm:resolution1}
for $\Bu= \OmE(\Ef_*)$,
\item $\la_{T}$ : $\ETp* \to \OmT$, the inverse of the map $\cH'$ (\ref{cH'}),
\item $\la_{D}$ : $\EDp* \to \OmD$, the map
$\mu'$ (\ref{map-mu1}) and
\item $\la_{C}$ : $\ECp* \to  \OmJ $, the composition
$\ga^{-1} \circ  \vr$ of the inverse of the map
$\ga$ (\ref{map-ga}) with the map $\vr$ (\ref{vr}).
\end{itemize}

The results of the previous section can be represented in the
form of the following commutative diagrams of sheaves of DGLAs,
their modules, and morphisms
\begin{equation}
\begin{array}{ccc}
(\ETp*,0,[,]_{SN}) &\stackrel{\la_{T}}{\brarrow} &(\OmT, D, [,]_{SN})\\[0.3cm]
\downarrow^{\Elie}_{\,mod}  & ~  &     \downarrow^{L}_{\,mod} \\[0.3cm]
(\EA*,0)  &\stackrel{\la_{A}}{\bbrarrow} & (\OmA, D),\\[1cm]
(\OmD, D+\pa, [,]_{G}) &\stackrel{\la_{D}}{\blarrow} & (\EDp*,\pa, [,]_{G})\\[0.3cm]
\downarrow^{\cR}_{\,mod}  & ~  &     \downarrow^{\ER}_{\,mod} \\[0.3cm]
(\OmJ, D+\db) &\stackrel{\la_{C}}{\bblarrow} & (\ECp*,\db),
\end{array}
\label{diag-T-D}
\end{equation}
where the horizontal arrows correspond to embeddings
of the sheaves of DGLAs (resp. of DGLA modules) constructed in
the previous section. These embeddings are quasi-isomorphisms
by theorems \ref{thm:resolution1},
\ref{thm-ona}  and propositions
\ref{resol-T-A}, \ref{mu-EDp}\,.

Next, due to claims {\it 1} and {\it 2} in theorem \ref{thm:kontsevich} we have
a fiberwise $L_\infty$-quasi-iso\-mor\-phism
\begin{equation}\label{cal-K}
\cK :  (\OmT,0,[,]_{SN}) \brarrow (\OmD,\pa,[,]_G)
\end{equation}
from the sheaf of DGLAs $(\OmT,0,[,]_{SN})$ to
the sheaf of DGLAs $(\OmD,\pa,[,]_G)$\,.
Composing $\Linf$-quasi-iso\-mor\-phism (\ref{cal-K}) with the
action of $\OmD$ on $\OmJ$ we get an $\Linf$-module structure
on $\OmJ$ over $\OmT$.

Due to claims {\it 1} and {\it 2} in theorem \ref{thm:shoikhet} we have
a fiberwise $\Linf$-quasi-iso\-mor\-phism
\begin{equation}\label{cal-S}
\cS :  (\OmJ,\db) \bbrarrow (\OmA,0)
\end{equation}
from the sheaf of $\Linf$-modules $\OmJ$ to the sheaf of
DGLA modules $\OmA$ over $\OmT$\,.

Thus we get the following commutative diagram
\begin{equation}
\begin{array}{ccc}
(\OmT, 0, [,]_{SN}) &\stackrel{\cK}{\brarrow} &(\OmD, \pa, [,]_{G})\\[0.3cm]
\downarrow^{L}_{\,mod}  & ~  &     \downarrow^{\cR}_{\,mod} \\[0.3cm]
(\OmA, 0)  &\stackrel{\cS}{\bblarrow} & (\OmJ, \db),
\end{array}
\label{diag-K-Sh}
\end{equation}
where by commutativity we
mean that $\cS$ is a $\Linf$-morphism of the sheaves of $\Linf$-modules
$(\OmJ, \db)$ and $(\OmE, 0)$ over the sheaf of DGLAs
$(\OmT, 0, [,]_{SN})$
and the $\Linf$-module structure on $(\OmJ, \db)$
over $(\OmT, 0, [,]_{SN})$ is obtained by composing
the $\Linf$-morphism $\cK$ with the action $\cR$
(see $3.4$ in \cite{D2}) of
$(\OmD, \pa, [,]_{G})$ on $(\OmJ, \db)$\,.

Let us now restrict ourselves to an open subset
$V\subset M$ such that $E\BV$ is trivial. Over any such
subset the $E$-de Rham differential (\ref{E-dif})
is well defined for either of the sheaves
$\OmA$, $\OmT$, $\OmJ$, and $\OmD$\,. Furthermore,
since the $\Linf$-quasi-isomorphisms (\ref{cal-K}) and (\ref{cal-S})
are fiberwise we can add to all the differentials
in diagram (\ref{diag-K-Sh}) the $E$-de Rham differential
(\ref{E-dif}). Thus we get a new commutative diagram
\begin{equation}
\begin{array}{ccc}
(\OmT\BV, \Edif, [,]_{SN}) &\stackrel{\cK}{\brarrow} &
(\OmD\BV, \Edif+\pa, [,]_{G})\\[0.3cm]
\downarrow^{L}_{\,mod}  & ~  &     \downarrow^{\cR}_{\,mod} \\[0.3cm]
(\OmA\BV, \Edif)  &\stackrel{\cS}{\bblarrow} & (\OmJ\BV, \Edif + \db)
\end{array}
\label{diag-V}
\end{equation}
of the $\Linf$-morphism $\cK$ and the
morphism of $\Linf$-modules $\cS$.

We claim that
\begin{prop}
The $\Linf$-morphism $\cK$
and the morphism of $\Linf$-modules $\cS$
in (\ref{diag-V}) are quasi-isomorphisms.
\end{prop}
\begin{proof}
This statement follows easily from
the standard argument of the spectral sequence.
Indeed, we can naturally regard $\OmT$ and $\OmD$
(resp. $\OmJ$ and $\OmA$) as sheaves of
double complexes and the exterior degree
provides us with the following descending filtration
$$
F^p(\OmE^{*}(\Bu)) = \bigoplus_{k\ge p} \OmE^k(\Bu)\,,
$$
where $\Bu$ is either $\Tp$ or $\Dp$ (resp.
$\Jp$ or $\Ef$).

The corresponding versions of Vey's \cite{V} and
Hoch\-schild-Kos\-tant-Ro\-sen\-berg-Connes-Te\-le\-man
 \cite{Co}, \cite{HKR}, \cite{Tel} theorems for
$\R^d_{formal}$ imply that $\cK$ (resp. $\cS$)
induces a quasi-isomorphism on the level of
$E_0$. Hence, $\cK$ (resp. $\cS$) induces a
quasi-isomorphism on the level of $E_{\infty}$.
The standard snake lemma argument of homological
algebra implies that $\cK$ (resp. $\cS$) in (\ref{diag-V}) is
a quasi-isomorphism.
\end{proof}

On the open subset $V$ we can represent the Fedosov
differential (\ref{DDD}) in the following
(non-covariant) form
\begin{equation}
\label{d+B}
D= \Edif + B \,\cdot \,,
\end{equation}
$$
B=\sum^{\infty}_{p=0} \xi^i B^k_{i;j_1 \dots j_p}(x) y^{j_1} \dots
y^{j_p} \frac{\pa}{\pa y^k}\,.
$$
If we regard $B$ as a section in $\OmE^1(V,\Tp^0)$ then
the nilpotency condition $D^2=0$ says that $B$ is a
Maurer-Cartan section of the sheaf of DGLAs
$(\OmE(\Tp)\BV, \Edif,[,]_{SN})$\,.
In the terminology of section $2$ in \cite{D2}
this means that the sheaf of DGLAs
$(\OmE(\Tp)\BV, $ $D ,[,]_{SN})$ is obtained
from
$(\OmE(\Tp)\BV, \Edif, [,]_{SN})$ via the twisting
procedure by the Maurer-Cartan element $B$\,.

According to proposition $1$ in section $2$
of \cite{D2} the element
$$
B_D=\sum_{k=1}^{\infty}\frac{1}{k!} \cK_k(B, \dots, B)
$$
is a Maurer-Cartan section of $(\OmD\BV, \Edif + \pa, [,]_{G})$\,.
Moreover, due to claim {\it 3} in theorem \ref{thm:kontsevich}
$$
B_D = B\,,
$$
where $B$ is viewed as a
section of the sheaf $\OmE^1(\Dp^0)\BV$\,.

Thus twisting the $\Linf$-quasi-isomorphism $\cK$
in (\ref{diag-V})
by the Maurer-Cartan element $B$ we get
the $\Linf$-quasi-isomorphism
$$
\cK^{tw} \,:\,(\OmT\BV, D ,[,]_{SN}) \brarrow
(\OmD\BV , D+\pa, [,]_{G})\,.
$$

Since the DGLA module structure on  $\OmA$ over
$\OmT$ (resp. on $\OmJ$ over $\OmD$) is honest
the twist by the Maurer-Cartan element
described in section $2$ of \cite{D2} do
not change these structures.
Hence, by virtue of propositions $3$ and
$4$ in \cite{D2} the twisting procedure turns diagram
(\ref{diag-V}) into the commutative diagram
\begin{equation}
\begin{array}{ccc}
(\OmT\BV, D, [,]_{SN}) &\stackrel{\tcK}{\brarrow} &
(\OmD\BV, D+\pa, [,]_{G})\\[0.3cm]
\downarrow^{L}_{\,mod}  & ~  &     \downarrow^{R}_{\,mod} \\[0.3cm]
(\OmA\BV, D)  &\stackrel{\tcS}{\bblarrow} & (\OmJ\BV, D + \db),
\end{array}
\label{diag-V1}
\end{equation}
where $\tcS$ is a $\Linf$-quasi-isomorphism obtained from $\cS$
by twisting via the Maurer-Cartan section $B$ of
the sheaf of DGLAs $(\OmT\BV, \Edif, [,]_{SN}) $\,.

We claim that the $\Linf$-morphism $\tcK$ (resp. $\tcS$)
does not depend on the choice of the
trivialization of $E$ over $V$ and hence is
a well-defined $\Linf$-morphism of sheaves of
DGLAs (resp. sheaves of DGLA modules).
Indeed, the term in $B$ that depends on the
choice of the trivialization of $E$
is linear in the fiber coordinates $y^i$\,.
But due to claim $4$ in theorem \ref{thm:kontsevich}
and claim $3$ in theorem \ref{thm:shoikhet} this term
contribute neither to $\tcK$ nor to $\tcS$\,.

Thus the $\Linf$-quasi-isomorphisms $\tcK$ and $\tcS$ are well defined
and we arrive at the following commutative diagram
\begin{equation}
\begin{array}{ccc}
(\OmT, D, [,]_{SN}) &\stackrel{\tcK}{\brarrow} &
(\OmD, D+\pa, [,]_{G})\\[0.3cm]
\downarrow^{L}_{\,mod}  & ~  &     \downarrow^{\cR}_{\,mod} \\[0.3cm]
(\OmA, D)  &\stackrel{\tcS}{\bblarrow} & (\OmJ, D + \db).
\end{array}
\label{diag-M}
\end{equation}

Assembling diagrams (\ref{diag-T-D}) and (\ref{diag-M}) we
get the desired chain (\ref{diag-thm}) of quasi-isomorphisms between
the sheaves of DGLA modules
$(\ETp*, \EA*)$ and $(\EDp*,$ $ \ECp*)$\,.
It is obvious from the construction that
the terms and the quasi-isomor\-phisms
of the resulting diagram (\ref{diag-thm})
are functorial in the pair $(E,\pa^E)$,
where $\pa^E$ is a torsion-free
connection on $E$.
Thus, theorem \ref{thm:formality} is proved. $\Box$

\subsection{Applications of the formality theorem}

The obvious applications of the formality theorem
for $E$-chains are related to the deformations associated with
Poisson Lie algebroids. Namely, theorem \ref{thm:formality}
allows us to get an elegant description
of the Hochschild homology and the traces
of these deformations.

First, we recall that
\begin{dfn}
A Lie algebroid $(E,M, \rho)$ equipped
with an $E$-bivector $\pi \in \G(M,\ETp{1})$
satisfying the Jacobi identity
\begin{equation}
\label{Jacobi}
[\pi,\pi]_{SN}=0
\end{equation}
is called a Poisson Lie algebroid.
\end{dfn}
Following \cite{NT} a quantization of
a Poisson Lie algebroid is a construction
of an element
\begin{equation}\label{Pi}
\Pi \in \G(M, \EDp{1})[[\h]]
\end{equation}
satisfying the condition of the
classical limit
\begin{equation}
\label{cl-lim}
\Pi = 1\otimes 1 ~{\rm mod}~ \h\,, \qquad
\Pi - t (\Pi) = \h \pi ~ {\rm mod}~\h \,,
\end{equation}
and the ``associativity'' condition
\begin{equation}\label{assoc}
[\Pi, \Pi]_G = 0\,.
\end{equation}
Here $\h$ is an auxiliary variable and
$t$ denotes the (cyclic) permutation
of components of
$\Pi\in \G(M, \EDp{1})[[\h]]$ $=$ $\G(M, \Op\otimes \Op)[[\h]]$\,.

Furthermore, two deformations $\Pi$ and
$\Pi'$ of $(E,M,\rho, \pi)$ are called \emph{equivalent}
if there exists a formal power series
$$
\Psi= 1 + \h\Psi_1 +\h^2 \Psi_2 + \, \dots\, \in \G(M,\Op)[[\h]]
$$
such that
\begin{equation}\label{equiv}
(\D \Psi)\,\Pi' = \Pi\,(\Psi\otimes \Psi)\,,
\end{equation}
where $\D$ is the coproduct
(\ref{coprod}) in $\Op$.

Thanks to the formality theorem
for the sheaf of DGLAs $\EDp*$ proved
in \cite{C} we have a bijective
correspondence between the moduli
spaces of Maurer-Cartan elements of
the DGLA $\hbar\G(M, \ETp*)[[\h]]$ of $E$-polyvector fields
and the DGLA $\hbar\G(M, \EDp*)[[\h]]$ of $E$-polydifferential
operators. 
In other words, if we consider the cone
\begin{equation}\label{cone}
\begin{array}{c}
\pi_{\h} = \h \pi + \h^2 \pi_1 +
+\h^3 \pi_2 + \dots\,, \\[0.3cm]
[\pi_{\h}, \pi_{\h}]_{SN}=0\,, \\[0.3cm]
\pi_i\in \G(M, \ETp{1})
\end{array}
\end{equation}
of formal power series in
$\h$ acted upon by the Lie
algebra $\h \G(M, E)[[\h]]$
\begin{equation}
\label{action}
\pi_{\h} \to [u, \pi_{\h}]\,, \qquad
u \in \h\G(M, E)[[\h]]\,,
\end{equation}
then
\begin{cor}
The deformations
(\ref{Pi}) associated with a Poisson Lie algebroid
$(E,M, \rho, \pi)$ modulo the relation
(\ref{equiv})
are in a bijective correspondence
with the points of the cone (\ref{cone})
modulo the action (\ref{action}) of
the prounipotent group corresponding to the
Lie algebra $\h\G(M, E)[[\h]]$\,. $\Box$
\end{cor}
An orbit $[\pi_{\h}]$ on the cone (\ref{cone})
corresponding to a deformation $\Pi$ (\ref{Pi})
is called \emph{the class of the deformation} and
any point $\pi_{\h}$ of this orbit is called
\emph{a representative} of the class.

Given a deformation $\Pi$ (\ref{Pi}) associated with
a Poisson Lie algebroid $(E,M, \rho, \pi)$ one
can define Hochschild chain complex of
this deformation as the graded vector space
\begin{equation}
\label{Hoch-hom}
\G(M, \ECp*)[[\h]]
\end{equation}
equipped with the differential
$$
\ER_{\Pi} : \ECp* \to \ECp{*+1}\,.
$$
Furthermore, one defines the Hochschild
cochain complex of the deformation $\Pi$
as the graded vector space
\begin{equation}
\label{Hoch-cohom}
\G(M, \EDp*)[[\h]]
\end{equation}
equipped with the differential
$$
[\Pi,\, ] : \EDp* \to \EDp{*+1}\,.
$$

Due to claim $5$ of proposition $2$ in \cite{D2},
claim $5$ of proposition $3$ in \cite{D2}, and
theorem \ref{thm:formality}
we get the following result:
\begin{cor}
\label{Hoch-comp}
Let $\Pi$ be a deformation associated
with a Poisson Lie algebroid $(E,M,\rho, \pi)$
and let $\pi_{\h}$ be a representative of
the class of this deformation. Then
the complex of Hochschild cohomology
(\ref{Hoch-cohom}) of the deformation
$\Pi$ is quasi-isomorphic to the complex
of $E$-polyvector fields
\begin{equation}
\label{E-fields}
(\G(M, \ETp*)[[\h]], [\pi_{\h},\, ])
\end{equation}
with the differential $[\pi_{\h},\,]$\,.
The complex of Hochschild homology
(\ref{Hoch-hom}) of the deformation
$\Pi$ is quasi-isomorphic to the complex of
$E$-forms
\begin{equation}
\label{Forms}
(\OmE(M)[[\h]], \Elie_{\pi_{\h}})
\end{equation}
with the differential $\Elie_{\pi_{\h}}$. $\Box$
\end{cor}

Given a deformation $\Pi$ (\ref{Pi}) associated with
a Poisson Lie algebroid $(E,M, \rho, \pi)$ one
can define \emph{a trace}\footnote{This notion 
is very important for formulations 
of various versions of algebraic index theorems 
for deformations associated with Poisson Lie 
algebroids \cite{NT}.} 
of the deformation
$\Pi$ as an $\R[[\h]]$-linear functional
\begin{equation}
\label{trace}
tr: \mO(M)[[\h]] \to \R[[\h]]
\end{equation}
satisfying the following condition
\begin{equation}
\label{trace1}
tr(j(\Pi) - j (t(\Pi))) =0\,,
\qquad \forall~j \in \G(M,\EJp{1}) \cap \ker \n^G\,.
\end{equation}

It is not hard to see that
corollary \ref{Hoch-comp} implies the
following statement:
\begin{cor}
Let $\Pi$ be a deformation associated
with a Poisson Lie algebroid $(E,M,\rho, \pi)$
and let $\pi_{\h}$ be a representative of
the class of this deformation. Then
the vector space of traces of the deformation $\Pi$ is isomorphic to
the vector space of continuous $\R[[\h]]$-linear
$\R[[\h]]$-valued functionals on $\mO(M)[[\h]]$
vanishing on all functions $f\in \mO(M)[[\h]]$ of the
following form
$$
f = j(\pi_{\h})\,, \qquad j \in \G(M,\EJp{1}) \cap \ker \n^G\,,
$$
where $\pi_{\h}$ is viewed as a series $E$-bidifferential
operators. $\Box$
\end{cor}

\section{Formality theorems for holomorphic Lie algebroids}\label{sect:holo}

Let now $M$ be a complex manifold. Let us write $T_\C M=T^{1,0}\oplus T^{0,1}$ for the
decomposition of the (complexified) tangent bundle as the sum of the holomorphic tangent bundle
and anti-holomorphic tangent bundle.
We denote by $\mO_M$ the structure sheaf of holomorphic functions on $M$ and
by $z^{\al}$ local coordinates on $M$\,.
We have to adapt the definition of holomorphic Lie algebroids:
\begin{dfn}
A holomorphic Lie algebroid over a complex manifold $M$
is a holomorphic vector bundle $E$ of finite rank
whose sheaf of sections is a sheaf of Lie algebras
equipped with a $\mO_M$-linear map of sheaves
of (complex) Lie algebras
$$
\rho:E\to T^{1,0}\,,
$$
satisfying the same condition described
(for the smooth case) in formula (\ref{eq:la}).
\end{dfn}
\noindent{\bf Remark.}
This notion is different from the one of a complex Lie algebroid that we introduced in the remark of 
subsection \ref{subsec:main}. 

\subsection{Algebraic structures and the main theorem}

Let $E$ be a holomorphic Lie algebroid.
As in section 1, one can define the
following sheaves (which are also holomorphic vector
bundles):
\begin{itemize}
\item $\ETp*$ is the sheaf of $E$-polyvector fields.
We regard $\ETp*$ as a sheaf of DGLAs
with the vanishing differential and
with the Lie bracket $[,]_{SN}$ defined as in
(\ref{S-N}), (\ref{eq:SN}).

\item $\EA*$ is the sheaf of $E$-differential
forms with converted grading:
\begin{equation}\label{EA1}
\EA* = \wedge^{-*} E^{\vee}\,, \qquad\qquad  \EA{0} =\mO_M\,.
\end{equation}
We regard $\EA*$ as a sheaf of DGLA modules
over $\ETp*$ with the vanishing differential
and with the action $\Elie$ defined as
in (\ref{C-Weil}).

\item $\EDp*$ is a sheaf of $E$-polydifferential
operators. We regard $\EDp*$ as a sheaf of DGLAs
with the bracket $[,]_G$ and the differential
$\pa$ defined as in (\ref{Gerst}) and
(\ref{pa}). Notice that the tensor product
of sections (over $\mO_M$) of $\EDp*$
turns $\EDp*$ into a sheaf of DGAs.

\item $\EJp*$ is the sheaf of $E$-polyjets
$$
\EJp{~} = \bigoplus_{k\ge 0} \EJp{k}\,, \quad \qquad
\EJp* :=
Hom_{\mO_M}(\EDp*, \mO_M)\,,
$$
which we regard as a sheaf
of DGLA modules over $\EDp*$ with the
action $\ES$ and the differential $\db$
defined as in (\ref{ono}) and (\ref{+db})\,.
The sheaf $\EJp*$ is also equipped
with the Grothendieck connection
\begin{equation}\label{eq:gro1}
\n^G : T^{1,0} \otimes \EJp* \mapsto \EJp*\,,
\qquad
\n^G_u(j)(P):=\rho(u)(j(P))-j(u \bul P)\,,
\end{equation}
where $u\in \G(U,T^{1,0})$ is a holomorphic
vector field,  $P\in \G(U,\EDp{k})$,
$j\in \G(U,\EJp{k})$ and the operation $\bul$
is defined in (\ref{bullet})\,.
The connection (\ref{eq:gro1}) is
compatible the DGLA module
structure on $\EJp*$\,.

\item $\ECp*$ is the graded sheaf of $\n^G$-flat
$E$-polyjets with converted grading
\begin{equation}\label{ECp1}
\ECp* := \ker \n^G\cap \EJp{-*}\,.
\end{equation}
Due to the compatibility of
the Grothendieck connection (\ref{eq:gro1})
with the DGLA module structure on $E$-polyjets
$\ECp*$ can be viewed as a sheaf of
DG modules over sheaf of DGLAs
$\EDp*$\,. We refer  to $\ECp*$ as
a sheaf of Hochschild $E$-chains or
$E$-chains for short.
\end{itemize}

The main result of this section can be formulated
as follows:
\begin{thm}\label{ya-Gilles}
For any holomorphic Lie algebroid
$E$ over a complex manifold $M$
the sheaves of DGLA modules
$(\ETp*, \EA*)$ and $(\EDp*, \ECp*)$
are quasi-isomorphic.
\end{thm}

Omitting the sheaves of DGLA modules $\EA*$ and
$\ECp*$ in the above theorem we get
the following corollary:
\begin{cor}\label{ya-Gilles1}
For any holomorphic Lie algebroid
$E$ over a complex manifold $M$
the sheaves of DGLAs
$\ETp*$ and $\EDp*$
are quasi-isomorphic.
\end{cor}
We would like to mention that
this corollary is parallel to the result of
A. Yekutieli \cite{Ye}, who
proved this statement for the
 tangent Lie algebroid
$TM \to M$ of any smooth algebraic
variety over a field $\K$ for which
$\R\subset \K$\,.

Notice that applying theorem \ref{ya-Gilles}
to the tangent algebroid $T^{1,0}M \to M$
we prove the following version of
Tsygan's formality conjecture for
complex manifolds:
\begin{thm}\label{Gilles1}
For any complex manifold $M$
the sheaf of DGLA modules $C^{poly}(M)$
of Hochschild chains over the
sheaf $D_{poly}(M)$ of (holomorphic)
polydifferential operators
is formal. $\Box$
\end{thm}

The proof of theorem \ref{ya-Gilles} occupies the
rest of the section.

\subsection{Fedosov resolutions}

First, we observe that any holomorphic Lie algebroid
$E$ can be viewed as a complex Lie algebroid in the sense
of the remark in subsection \ref{subsec:main},
where the anchor map is naturally extended to
the smooth sections of $E$. 
It is clear that the sheaf of
Lie algebras $T^{0,1}$ acts on $E$ and that this
action commutes with $\rho$ as $\rho$ is holomorphic.
Thus we get
\begin{prop}
Let $F$ be the smooth vector bundle $F=E \oplus T^{0,1}$.
Then $F$ is a complex Lie algebroid over $M$ with the anchor map $\rho_F$ :
$F \to T^{1,0} \oplus T^{0,1}$ given by $\rho_F\BE=\rho$ and
$\rho_F\BTO = \id\, :  T^{0,1} \to T^{0,1}$. $\Box$
\end{prop}

For a holomorphic
vector bundle $\Bu$ over $M$
we consider the sheaf of smooth
$F$-differential forms with values in $\Bu$:
\begin{equation}\label{OmF-Bu}
\OmEF(\Bu)=\bigoplus_{p,q}\OmEF^{p,q}(\Bu)\,,
\end{equation}
$$
\OmEF^{p,q}(\Bu) =
\wedge^p E^{\vee}\otimes \wedge^q T^{*0,1}M\otimes \Bu
$$
For sections
\begin{equation}\label{sect-OmF-Bu}
a=\sum_{p,q} a_{i_1\dots i_p\, ;\, \al_1, \dots, \al_q}(z,\bz)
\xi^{i_1} \dots \xi^{i_p} d \bz^{\al_1} \dots d \bz^{\al_q}\,,
\end{equation}
$$
a_{i_1\dots i_p\, ;\, \al_1, \dots, \al_q}(z,\bz) \in
\G^{\,\rm smooth}(U,\Bu)
$$
of $\OmEF(\Bu)$ we reserve
the local basis $\{\xi^i\}$ of
anti-commuting fiber coordinates
on $E$ and the local basis
$\{ d \bz^{\al} \}$ of antiholomorphic
exterior forms on $M$\,. We denote by
$\brd$ the Dolbeault differential
\begin{equation}\label{Dolbo}
\brd = d \bz^{\al}\pa_{\bz^{\al}} :
\OmEF^{p,*}(\Bu) \mapsto \OmEF^{p, *+1}(\Bu)\,.
\end{equation}

It is obvious that the
(DG) algebraic structures on the sheaves
$\ETp*$, $\EA*$, $\EDp*$, and $\EJp*$,
can be naturally extended to
the sheaves $\OmEF^{0,*}(\ETp*)$,
$\OmEF^{0,*}(\EA*)$, $\OmEF^{0,*}(\EDp*)$,
and $\OmEF^{0,*}(\EJp*)$\,. Similarly,
the Grothendieck
connection (\ref{eq:gro1}) on $\EJp*$
extends to the operator
\begin{equation}\label{eq:gro11}
\n^G : T^{1,0} \otimes \OmEF^{0,*}(\EJp*)
\mapsto \OmEF^{0,*}(\EJp*) \,,
\end{equation}
which is compatible with the
action $\ES$ of $\OmEF^{0,*}(\EDp*)$
on $\OmEF^{0,*}(\EJp*)$ and with
the differential $\db$ on  $\OmEF^{0,*}(\EJp*)$\,.

Since $\ETp*$, $\EA*$, $\EDp*$, and $\EJp*$
are holomorphic vector bundles it makes
sense to speak about
the Dolbeault differential (\ref{Dolbo})
\begin{equation}\label{Dol}
\brd : \OmEF^{0,*}(\Bu) \mapsto \OmEF^{0, *+1}(\Bu)\,,
\end{equation}
for $\Bu$ being either $\ETp*$, $\EA*$, $\EDp*$,
or $\EJp*$\,. It is obvious that
$\brd$ is compatible with the
(DG) algebraic structures on $\OmEF^{0,*}(\Bu)$
and with the Grothendieck connection
(\ref{eq:gro11}) on $\OmEF^{0,*}(\EJp*)$\,.

Furthermore, due to the $\brd$-Poincar\'e
lemma we have
\begin{prop}\label{Gilles}
If $\Bu$ is either $\ETp*$, $\EA*$,
$\EDp*$, or $\EJp*$ then the
canonical inclusion of sheaves
\begin{equation}\label{Dol1}
\inc : \Bu \hookrightarrow \OmEF^{0,*}(\Bu)
\end{equation}
is a quasi-isomorphism of complexes
of sheaves $(\Bu,0)$ and
$(\OmEF^{0,*}(\Bu), \brd )$\,.
The inclusion $\inc$ is compatible with the
(DG) algebraic structures on $\Bu$,
and $\OmEF^{0,*}(\Bu)$, and
with the Grothendieck connection
(\ref{eq:gro1}), (\ref{eq:gro11}). $\Box$
\end{prop}
Due to this proposition it suffices to
prove that the sheaves of DGLA modules
$(\OmEF^{0,*}(\ETp*),$ $\OmEF^{0,*}(\EA*))$,
and $(\OmEF^{0,*}(\EDp*),$ $\OmEF^{0,*}(\EJp*))$
are quasi-isomorphic.
To show this we follow the lines
of section 2 and introduce the formally completed
symmetric algebra $\SE$ of the dual bundle $E^{\vee}$ and
(holomorphic) bundles $\Tp,$ $\Dp,$ $\Ef,$ $\Jp$ associated with $\SE$
(see page \pageref{sect-Tp}). As in section 2,
$\Tp$ and $\Dp$ are sheaves of DGLAs while
$\Ef$ and $\Jp$ are sheaves of DGLA modules
over $\Tp$ and $\Dp$, respectively. $\Dp$ is
also a sheaf of DGA's.

Next, we consider sheaves of smooth
$F$-differential forms with values
in the bundles $\SE$ $\Tp$, $\Dp$, $\Ef$,
and $\Jp$\,. It is clear that
the sheaves $\OmSF$, $\OmAF$, $\OmTF$, $\OmDF$,
and $\OmJF$ acquire the corresponding
(DG) algebraic structures and the
Dolbeault differential (\ref{Dolbo})
is obviously compatible with these
structures.

Furthermore, we have the following
obvious analogue of proposition \ref{ono1}
\begin{prop}\label{ono11}
The sheaf $\OmEF(\Tp^0)$ of $F$-forms with values
in fiberwise vector fields is a sheaf of graded Lie
algebras. The sheaves
$\OmSF$, $\OmAF$, $\OmTF$, $\OmDF$, and $\OmJF$
are sheaves of modules over $\OmEF(\Tp^0)$ and
the action of $\OmEF(\Tp^0)$ is
compatible with the DG algebraic structures
on $\OmSF$, $\OmAF$, $\OmTF$, $\OmDF$, $\OmJF$
and with the Dolbeault differential
(\ref{Dolbo}). $\Box$
\end{prop}
Due to this proposition one can extend the following
differential
$$
\delta := \xi^i\frac{\pa}{\pa y^i}:
\OmEF^{*,q}(\SE) \to\OmEF^{*+1,q}(\SE)
$$
of the sheaf of algebras $\OmSF$ to
the sheaves $\OmTF, \OmDF,$ $\OmAF$ and $\OmJF$ so
that $\de$ is compatible with the (DG) algebraic
structures on $\OmTF$, $\OmAF$,
$\OmDF$, and $\OmJF$, and with the differential
$\brd$ (\ref{Dolbo}).
Here $\{ y^i \}$ (resp. $\{ \xi^i \}$)
denote commuting (resp. anticommuting) fiber coordinates
of the bundle $E$\,.

\smallskip

We now have an analogue of proposition \ref{thm:delta}
\begin{prop}
\label{thm:delta-comp}
For $\Bu$ being either of the sheaves
$\SE$, $\Ef$, $\Tp$ or $\Dp$ and $q \geq 0$,
$$
H^{\ge 1}(\OmEF^{*,q}(\Bu), \de) =0\,.
$$
Furthermore,
\begin{equation}
\begin{array}{c}
H^0(\OmEF^{*,q}(\SE), \de) \cong \OmEF^{0,q}(M,\mO_M)\,, \\[0.3cm]
H^0(\OmEF^{*,q}(\Ef_*), \de) \cong  \OmEF^{0,q}(M,\EA*)\,, \\[0.3cm]
H^0(\OmEF^{*,q}(\Tp^*), \de) \cong  \OmEF^{0,q}(M,\wedge^{*+1}(E))
\end{array}
\end{equation}
as sheaves of (graded) commutative algebras
and
\begin{equation}
\qquad
H^0(\OmEF^{*,q}(\Dp^*), \de) \cong  \OmEF^{0,q}(M,\otimes^{*+1}(S(E)))
\end{equation}
as sheaves of DGAs
over $\mO_M$.
\end{prop}
\begin{proof}
As in proposition \ref{thm:delta}
is suffices to construct an operator
($q\geq 0$)
\begin{equation}
\ka: \OmEF^{*,q}(\Bu) \to \OmEF^{*-1,q}(\Bu)
\label{kappa-comp}
\end{equation}
such that for any section $u$
of $\OmEF(\Bu)$ equation
\begin{equation}
\label{eq:homotopy1}
u=\delta\kappa (u)+\kappa\delta (u) + \cH(u)\,,
\end{equation}
is still true, where now
\begin{equation}
\label{cH1}
\cH(u) = u \Big|_{y^i=\xi^i=0}\,.
\end{equation}
and $y^i$ are as above
fiber coordinates on $E$\,.
As in the proof of proposition \ref{thm:delta}
we define $\ka$ on $\OmEF(\SE)$ by
equation (\ref{kappa}) and then extend it to
$\OmTF$, $\OmAF$, and $\OmDF$ in the
componentwise manner.
\end{proof}

Let us choose a connection $\pa^E$ on $E$
which is compatible with the complex
structure on $E$
\begin{equation}\label{pa-E}
\pa^E = \Edif + \brd + \xi^{i} \G_i :
\OmEF^{*}(E) \to \OmEF^{*+1}(E)\,,
\end{equation}
where $\xi^i \G_i$ is locally a section
of the sheaf $\OmE^1(End(E))$ and
$\Edif : \OmEF^{*,q}_M  \to \OmEF^{*+1,q}_M$
is defined in (\ref{E-dif})\,.

It is not hard to show that such a connection
always exists, and moreover, one can always choose
$\pa^E$ to be torsion free.

As in (\ref{nabla})
we extend (\ref{pa-E}) to a
derivation of
the DG sheaves $\OmSF$, $\OmAF$, $\OmTF$,
$\OmDF$, and $\OmJF$:
\begin{equation}\label{nabla1}
\n = \Edif + \G\cdot + \brd :\OmEF^{*}(\Bu) \to \OmEF^{*+1}(\Bu)\,,
\end{equation}
where $\Bu$ is either of the sheaves $\SE$, $\Ef$,
$\Tp$, $\Dp$, or $\Jp$,
$\G = -\xi^i \Gamma_{ij}^k y^j\frac{\pa}{\pa y^k}$\,,
$\Gamma_{ij}^k(x)$ are
Christoffel's symbols of the connection $\pa^E$
(\ref{pa-E}) and
$\G\cdot$ denotes the action of $\G$ on the
sections of the sheaves $\OmEF(\Bu)$\,.
Due to proposition \ref{ono11} $\n$ (\ref{nabla1}) is
compatible with the DG algebraic
structures on $\OmSF$, $\OmTF$, $\OmAF$,
$\OmDF$, and $\OmJF$, and since $\n$ is torsion
free
\begin{equation}\label{nado}
\n \de +\de \n = 0\,.
\end{equation}

\smallskip

Regarding (\ref{nabla1}) as a connection
on $\Bu$ one can see that the curvature
of (\ref{nabla1}) has the
components of type $(2,0)$ and $(1,1)$
\begin{equation}\label{nabla-R}
\n^2 = R^{2,0}+ R^{1,1}\,, \qquad
R^{2,0} = ( \Edif + \G )^2, \qquad
R^{1,1} = \bar{d}\G\,.
\end{equation}

\smallskip

We now prove the existence of a complex Fedosov differential $D$:
\begin{thm}
Let $\Bu$ be either of the sheaves
$\SE$, $\Ef$, $\Tp$, $\Dp$, or $\Jp$\,.
There exists a section
\begin{equation}
\label{eq:Acomp}
A=\sum_{s=2}^\infty\xi^k
A_{k,i_1\dots i_s}^j(z,\bz) y^{i_1}\cdots y^{i_s}\frac{\pa}{\pa y^j}
\end{equation}
of the sheaf $\OmEF^{1,0}(\Tp^0)$
and a section
\begin{equation}
\label{eq:Bcomp}
\bar{A}=\sum_{s=2}^\infty d \bar{z}^\alpha
\bar{A}_{\alpha,i_1\dots i_s}^j(z,\bz) y^{i_1}\cdots y^{i_s}\frac{\pa}{\pa y^j}
\end{equation}
of the sheaf $\OmEF^{0,1}(\Tp^0)$
such that the derivation
\begin{equation}
\label{DDDcomp}
D:= \n  - \delta  + A \cdot + \bar{A} \cdot : \OmEF^{*}(\Bu)
\to \OmEF^{*+1}(\Bu)
\end{equation}
is $2$-nilpotent ($D^2 = 0$)
and compatible with the
DG  algebraic structure on $\OmEF(\Bu)$\,.
\end{thm}
\begin{proof}
Let us rewrite $D=D^{1,0}+D^{0,1}$ with
$$
D^{1,0}=\Edif + \G\, \cdot\, - \delta + A \cdot\,,
\qquad
D^{0,1}=\bar{d}+ \bar{A} \cdot
$$
and try to mimic the proof of theorem \ref{G-Fed}.

Due to (\ref{nado}) and (\ref{nabla-R})
the condition $(D^{1,0})^2=0$ is
equivalent to the equation
$$
R^{2,0} +  (\Edif + \G \, \cdot \, ) A - \de A  + \frac12 [A,A]_{SN} =0\,.
$$
This equation has a solution obtained
by iterations of the
following equation (with respect to the degrees
in fiber coordinates $y_i$'s)
$$
A=\kappa R^{2,0} + \kappa( (\Edif + \G \, \cdot \,) A+\frac12[A, A]_{SN})
$$
(the proof is the same as for theorem \ref{G-Fed}).

Using (\ref{nabla-R}) once again we observe that
the condition $D^{1,0}D^{0,1}+D^{0,1}D^{1,0}=0$
is equivalent to
$$
R^{1,1}+\bar{d}A+ (\Edif + \G \, \cdot \, ) \bar{A} -\delta \bar{A} +[A,\bar{A}]_{SN} =0\,,
$$
which, using the same arguments, has a solution obtained by
iterations of the equation
$$
\bar{A}=\kappa(R^{1,1}+\bar{d}A+ (\Edif + \G \, \cdot \, ) \bar{A} +[A,\bar{A}]_{SN}).
$$
Indeed, denoting
$$
C^{1,1}=R^{1,1}+\bar{d}A+ (\Edif + \G \, \cdot \, ) \bar{A} -
\delta \bar{A} +[A,\bar{A}]_{SN}\,,
$$
and using that
$ \de A= R^{2,0} +  \Edif + \G \, \cdot \, A  + \frac12 [A,A]_{SN} $ ($(D^{1,0})^2=0$),
$\bar{d} R^{2,0}=0$ and $\delta R^{1,1}=0$
(Bianchi's identities for $\n$) we get
$$
(\Edif + \G \, \cdot \, ) C^{1,1}-\delta C^{1,1} + [A,C^{1,1}]=0.
$$
We have $\kappa C^{1,1}=0$ by construction of $\bA$ and so,
by the ``Hodge-de Rham'' decomposition (\ref{eq:homotopy1}), we have
$$
C^{1,1}= \kappa( (\Edif + \G \, \cdot \, ) C^{1,1} + [A, C^{1,1}]).
$$
The latter equation has the unique
vanishing solution, which gives the result.

Let us now check the condition $(D^{0,1})^2=0$.
This will be true if the section
$$
C^{0,2}=\bar{d} \bar{A} + \frac12[\bar{A}, \bar{A}] \in \OmEF^{0,2}(\Tp^0)
$$
is zero.
One has again $D^{1,0}C^{0,2}=0$ and $\kappa C^{0,2}=0$ because
it does not have $\xi$'s.
As before, one can conclude that $C^{0,2}=0$.

The compatibility of (\ref{DDDcomp}) with the corresponding
DG algebraic structures follows from proposition
\ref{ono11}.
\end{proof}

We now describe the cohomology
of the Fedosov differential $D$ for the
sheaves $\OmSF$, $\OmAF$, $\OmTF$, and $\OmDF$
\begin{thm}\label{cohom-comp}
Let $\Bu$ be either of the sheaves $\SE$, $\Ef$, $\Tp$, or $\Dp$
and $q \geq 0$. We have
$$H(\OmEF^{*}(\Bu), D) \cong
H(\OmEF^{0,*}(M,\Bu)\cap \ker \delta,\bar{d})\,.$$
as sheaves of (differential) graded (commutative) algebras.
\end{thm}
\begin{proof}
Let us consider the double complex
$(\OmEF^{*,*}(\Bu), D^{1,0}+D^{0,1})$.
Using the degree in the fiber coordinates $y^i$ we
introduce on this complex a decreasing filtration.
Applying the spectral sequence argument (as in the proof of
theorem \ref{thm:resolution1}) and using proposition
\ref{thm:delta-comp} we conclude that
for any $i\geq 0$, the cohomology of the complex
$(\OmEF^{*,i}(\Bu), D^{1,0})$ is concentrated in
degree $*=0$\,. Therefore,
\begin{equation}\label{Vasya}
H(\OmEF^{*}(\Bu), D) =
H(\OmEF^{0,*}(\Bu)\cap \ker D^{1,0}, D^{0,1})\,.
\end{equation}

Following the lines of the
proof of theorem \ref{thm:resolution1}
it is not hard to show that iterating the
equation
\begin{equation}\label{iter-u1}
\la(u) = u + \ka (\n \la(u) + A\cdot \la(u) +  \bA\cdot \la(u))\,,
\qquad u \in  \OmEF^{0,q}(U,\Bu) \cap \ker \de
\end{equation}
we get an isomorphism of sheaves
(of graded vector spaces)
\begin{equation}\label{lift1}
\la: \OmEF^{0,q}(\Bu ) \cap \ker \delta\to
\OmEF^{0,q}(\Bu)\cap \ker D^{1,0}\,,
\end{equation}
and moreover, the map $\la$ (\ref{lift1}) has
a natural inverse given by the map
$\cH$ (\ref{cH1}).

We claim that $\la$ gives a quasi-isomorphism
of complexes
$$
\la: (\OmEF^{0,*}(\Bu)\cap \ker \delta,\bar{d})\to
(\OmEF^{*}(\Bu), D).
$$
Indeed, due to (\ref{Vasya}) it suffices
to show that for any
$u \in \OmEF^{0,q}(U,\Bu)\cap \ker \delta$, one has
$$
\lambda (\bar{d}(u))=D^{0,1}\lambda(u).
$$
The term $\la (\bar{d}(u))$ is the only element in
$\OmEF^{0,q}(\Bu)$ such that $\cH(\la (\bar{d}(u)))=\bar{d}(u)$
and $D^{1,0} \la (\bar{d}(u)) =0$. It is clear that
$\cH(D^{0,1}\lambda(u))=\bar{d}(u)$
and one has
$$
D^{1,0} D^{0,1}\lambda(u)=-D^{0,1}D^{1,0} \lambda(u)=0,
$$
since map $\la$ (\ref{iter-u1}) lands
in $\ker D^{1,0}$\,.

\smallskip

The map $\la$ (\ref{lift1}) is compatible with
the corresponding multiplications in $\SE$, $\Ef$, $\Tp$, or
$\Dp$ since so is the map $\cH$ (\ref{cH1})\,.
The theorem is proved.
\end{proof}

\smallskip

It is not hard to prove the following
analogue of proposition \ref{resol-T-A} :
\begin{prop}\label{resol-T-A1}
The composition
\begin{equation}\label{cH'-comp}
\cH'=\nu \circ \cH :
 \OmEF^{0,*}(\Tp) \cap \ker D^{1,0} \to \OmEF^{0,*}(\ETp*)
\end{equation}
of the map
\begin{equation}\label{cHHH-comp}
\cH\,:\,  \OmEF^{0,*}(\Tp) \cap \ker D^{1,0}
\to  \OmEF^{0,*}(\Tp) \cap \ker \de
\end{equation}
with the identification
of the sheaves $\Tp^* \cap \ker \de$ and
$\ETp* \cong \wedge^{*+1} E$
\begin{equation}\label{nu-comp}
\nu : \Tp^* \cap \ker \de
\erarrow \ETp*
\end{equation}
is an isomorphism of the sheaves
of DGLAs
\begin{equation}\label{DGLA}
(\OmEF^{0,*}(\Tp) \cap \ker D^{1,0}, D^{0,1}, [,]_{SN}) \cong
(\OmEF^{0,*}(\ETp*), \brd, [,]_{SN})
\end{equation}
The map
\begin{equation}\label{cH-A1-comp}
\cH\,:\, \OmEF^{0,*}(\Ef_*) \cap \ker D^{1,0}
\to \OmEF^{0,*}(\EA*)
\end{equation}
is an isomorphism of the sheaves of
DGLA modules
$$
(\OmEF^{0,*}(\Ef_*) \cap \ker D^{1,0} , D^{0,1})  \cong
(\OmEF^{0,*}(\EA*), \brd)
$$
over the sheaf of DGLAs (\ref{DGLA})\,. $\Box$
\end{prop}
Thanks to equation (\ref{Vasya}) this
proposition implies that the map
$\cH'$ gives a quasi-isomorphism of
the sheaves of DGLAs
$(\OmEF^{*}(\Tp), D, [,]_{SN})$ and
$(\OmEF^{0,*}(\ETp*), \brd, [,]_{SN})$\,.

Playing with the PBW theorem for the Lie algebroids
(as we did in the proof of proposition \ref{mu=iso})
and with the cup product in the
sheaves $\Dp$ and $\EDp*$
(see equation (\ref{map-mu1})) one can
prove the following analogue of
proposition \ref{mu-EDp}
\begin{prop}\label{mu-EDp-comp}
The exists an isomorphism
of the sheaves of DGLAs
\begin{equation}\label{map-mu-comp}
\mu' : (\OmEF^{0,*}(\EDp*), \brd, [,]_G ) \erarrow
(\OmEF^{0,*}(\Dp) \cap \ker D^{1,0}, D^{0,1}, [,]_G)\,,
\end{equation}
which is compatible
with the DGA
structures on the sheaves $\OmEF^{0,*}(\EDp*)$ and
$\OmEF^{0,*}(\Dp)$. $\Box$
\end{prop}
Thanks to equation (\ref{Vasya}) this
proposition implies that the map
$\mu'$  (\ref{map-mu-comp}) gives a quasi-isomorphism of
the sheaves of DGLAs
$(\OmEF^{*}(\Dp), D, [,]_{G})$ and
$(\OmEF^{0,*}(\EDp*), \brd, [,]_{G})$\,.

Let us consider the map
\begin{equation}\label{map-ga-comp}
\ga : \OmEF^{0,q}(\Jp_*)
\to \OmEF^{0,q}(\EJp*)\,,
\qquad
\ga(j)(P)= (\mu' (P))(j)\Big|_{y^i=0}\,,
\end{equation}
where $ j\in \OmEF^{0,q}(U,\Jp_k)$ and
$P$ is a holomorphic section of $\EDp{k}$\,.

For this map
we have the following obvious
analogue of theorem \ref{thm-ona}
\begin{thm}\label{thm-ona-comp}
For any $q\ge 0$
\begin{equation}\label{Vasya1}
H^q(\OmEF^{*}(\Jp), D) =
H^q(\OmEF^{0,*}(\Jp)\cap \ker D^{1,0}, D^{0,1})\,.
\end{equation}
and the map $\ga$ (\ref{map-ga-comp}) provides
us with an isomorphism of
the sheaves of DGLA modules
\begin{equation}\label{H0D-J-comp}
\ga :
\OmEF^{0,*}(\Jp_*)
\erarrow
\OmEF^{0,*}(\EJp*)
\end{equation}
over the sheaf of DGLAs
$$
(\OmEF^{0,*}(\Dp) \cap \ker D^{1,0}, D^{0,1}, [,]_G)
\cong
(\OmEF^{0,*}(\EDp*), \brd, [,]_G )\,.
$$
The map $\ga$ sends the component $D^{1,0}$
to the Grothen\-dieck connection
(\ref{eq:gro11}) and
the component $D^{0,1}$ to the
Dolbeault differential $\brd$
(\ref{Dolbo}). $\Box$
\end{thm}

\subsection{End of the proof of theorem \ref{ya-Gilles}}

We have constructed the
following honest ($\Linf$-)quasi-isomorphisms of
the sheaves of DGLA modules
\begin{itemize}
\item $\la_{T}$ : $(\OmEF^{0,*}(M,\ETp*),\bar{d}, [,]_{SN}) \to
(\OmTF,D, [,]_{SN})$, \\[0.1cm]

\item $\la_{A}$ : $(\OmEF^{0,*}(M,\EA*),\bar{d})\to
(\OmEF(\Ef_*),D)$,\\[0.1cm]

\item $\la_{D}$ : $(\OmEF^{0,*}(M,\EDp*),\bar{d}, [,]_G) \to
(\OmDF,D, [,]_G)$, and \\[0.1cm]

\item $\la_{C}$ : $(\OmEF^{0,*}(M,\ECp*),\bar{d}) \to  (\OmJF,D) $.
\end{itemize}
Namely, the map $\la_T$ is the inverse of $\cH'$
(\ref{cH'-comp}) the map $\la_{A}$ is the inverse of
$\cH$ (\ref{cH-A1-comp}) $\la_D=\mu'$ (\ref{map-mu-comp}), and
$\la_C$ is composition of the identification
(\ref{ECp1}) and the inverse of $\ga$
(\ref{map-ga-comp}).

Our results can be summarized in
the following commutative diagrams
\begin{equation}
\label{ssylka}
\begin{array}{ccc}
(\OmEF^{0,q}(M,\ETp*),\bar{d},[,]_{SN}) &\stackrel{\la_{T}}{\brarrow} &(\OmTF, D, [,]_{SN})\\[0.3cm]
\downarrow^{\Elie}_{\,mod}  & ~  &     \downarrow^{L}_{\,mod} \\[0.3cm]
(\OmEF^{0,q}(M,\EA*),\bar{d})  &\stackrel{\la_{A}}{\bbrarrow} & (\OmAF, D),\\[1cm]
(\OmDF, D+\pa, [,]_{G}) &\stackrel{\la_{D}}{\blarrow} & (\OmEF^{0,q}(M,\EDp*),\bar{d}+\pa, [,]_{G})\\[0.3cm]
\downarrow^{\cR}_{\,mod}  & ~  &     \downarrow^{\ER}_{\,mod} \\[0.3cm]
(\OmJF, D+\db) &\stackrel{\la_{C}}{\bblarrow} &
(\OmEF^{0,q}(M,\ECp*),\bar{d}+\db),
\end{array}
\end{equation}
where the action $\ER$ is obtained from the action $\ES$
of $\OmEF^{0,*}(M,\EDp*)$ on $\OmEF^{0,q}(M,\EJp*)$ via
the identification (\ref{ECp1}).

Due to claims {\it 1} and {\it 2} in theorem \ref{thm:kontsevich}
and claims {\it 1} and {\it 2} in theorem \ref{thm:shoikhet}
we get the following commutative diagram
\begin{equation}
\begin{array}{ccc}
(\OmTF, 0, [,]_{SN}) &\stackrel{\cK}{\brarrow} &(\OmDF, \pa, [,]_{G})\\[0.3cm]
\downarrow^{L}_{\,mod}  & ~  &     \downarrow^{\cR}_{\,mod} \\[0.3cm]
(\OmAF, 0)  &\stackrel{\cS}{\bblarrow} & (\OmJF, \db),
\end{array}
\label{diag-K-Sh-holo}
\end{equation}
where by commutativity we
mean that $\cS$ is a morphism of the sheaves of $\Linf$-modules
$(\OmJF, \db)$ and $(\OmEF, 0)$ over the sheaf of DGLAs $(\OmTF, 0, [,]_{SN})$
and the $\Linf$-module structure on $(\OmJF, \db)$
over $(\OmTF, 0, [,]_{SN})$ is obtained by composing
the $\Linf$-isomorphism $\cK$ with the action $\cR$
(see $3.4$ in \cite{D2}) of
$(\OmDF, \pa, [,]_{G})$ on $(\OmJF, \db)$\,.

Let us now restrict ourselves to an open subset
$V\subset M$ such that $E\BV$ is trivial. Over any such
subset the $E$-de Rham differential (\ref{E-dif})
is well defined for either of the sheaves
$\OmAF$, $\OmTF$, $\OmJF$, and $\OmDF$\,. So again, we get a new commutative diagram
\begin{equation}
\begin{array}{ccc}
(\OmTF\BV, \Edif+\bar{d}, [,]_{SN}) &\stackrel{\cK}{\brarrow} &
(\OmDF\BV, \Edif+\bar{d}+\pa, [,]_{G})\\[0.3cm]
\downarrow^{L}_{\,mod}  & ~  &     \downarrow^{\cR}_{\,mod} \\[0.3cm]
(\OmAF\BV, \Edif+\bar{d})  &\stackrel{\cS}{\bblarrow} & (\OmJF\BV, \Edif+\bar{d}+ \db)
\end{array}
\label{diag-V-holo}
\end{equation}
in which the $\Linf$-morphism $\cK$ and the
morphism of $\Linf$-modules $\cS$ are
quasi-isomorphisms.

On the open subset $V$ we can represent the Fedosov
differential (\ref{DDD}) in the following
(non-covariant) form
\begin{equation}\label{d+B-holo}
D= \Edif +\bar{d}+ B\, \cdot\, + \bar{B}\,\cdot\,,
\end{equation}
$$
B=\sum^{\infty}_{p=0} \xi^i B^k_{i;j_1 \dots j_p}(z,\bz) y^{j_1} \dots
y^{j_p} \frac{\pa}{\pa y^k}\,,
$$
and
$$
\bar{B}=\sum^{\infty}_{p=0}d \bar{z}^\alpha
\bar{B}^k_{\alpha;j_1 \dots j_p}(z, \bz) y^{j_1} \dots
y^{j_p} \frac{\pa}{\pa y^k}\,,
$$
where the $z^\al$ are local coordinates on $M$\,.
If we regard $B+\bar{B}$ as a section of $\OmEF^1(\Tp^0)\BV$ then
the nilpotency condition $D^2=0$ says that $B+\bar{B}$ is a
Maurer-Cartan section of the sheaf of DGLAs
$(\OmEF(\Tp)\BV, \Edif+\bar{d},[,]_{SN})$\,.

Thus applying the twisting procedures developed
in section $2$ of \cite{D2} and using claim $3$
of theorem \ref{thm:kontsevich} we get
the following commutative diagram
\begin{equation}
\begin{array}{ccc}
(\OmT\BV, D, [,]_{SN}) &\stackrel{\tcK}{\brarrow} &
(\OmD\BV, D+\pa, [,]_{G})\\[0.3cm]
\downarrow^{L}_{\,mod}  & ~  &     \downarrow^{\cR}_{\,mod} \\[0.3cm]
(\OmA\BV, D)  &\stackrel{\tcS}{\bblarrow} & (\OmJ\BV, D + \db),
\end{array}
\label{diag-V1-comm}
\end{equation}
in which $\tcK$ is a $\Linf$-quasi-isomorphism of
the sheaves of DGLAs and $\tcS$ is a $\Linf$-quasi-isomorphism
of the sheaves of DGLA modules.

Due to claim $4$ in theorem \ref{thm:kontsevich}
and claim $3$ in theorem  \ref{thm:shoikhet}
the quasi-isomorphisms do not
depend on the trivialization of $E$
over $V$\,.

Thus we constructed the following commutative
diagram of sheaves of DGLAs, DGLA modules
and their $\Linf$-quasi-isomorphisms:
\begin{equation}
\begin{array}{ccc}
(\OmT, D, [,]_{SN}) &\stackrel{\tcK}{\brarrow} &
(\OmD, D+\pa, [,]_{G})\\[0.3cm]
\downarrow^{L}_{\,mod}  & ~  &     \downarrow^{\cR}_{\,mod} \\[0.3cm]
(\OmA, D)  &\stackrel{\tcS}{\bblarrow} & (\OmJ, D + \db),
\end{array}
\label{diag-V1-comp}
\end{equation}

Combining the diagrams in (\ref{ssylka}),
(\ref{diag-V1-comp}) together with the
proposition \ref{Gilles}
we see that the sheaves of DGLA modules
$(\ETp*, \EA*)$ and $(\EDp*, \ECp*)$
are connected by chain of quasi-isomorphisms.
Thus, theorem \ref{ya-Gilles} is proved. $\Box$

\section{Concluding remarks}

It would be interesting to prove
the corresponding version of the
algebraic index theorem \cite{NT}, \cite{TT1},
which should relate a cyclic chain in the
complex associated with a deformation
$\Pi$ (\ref{Pi}) to its principal part and
characteristic classes of the
Lie algebroid $(E, M, \rho)$.
It would be also interesting to investigate
how other characteristic classes
\cite{Cr}, \cite{Fer}, \cite{LSh} of Lie
algebroids could enter this picture.

Corollary \ref{ya-Gilles1} does not
in general give a chain of
quasi-isomorphisms between the
DGLAs $\G(M,\ETp*)$ and $\G(M,\EDp*)$
of global sections. However,
since the sheaves of smooth
forms $\OmEF^{0,*}(\ETp*)$ and
$\OmEF^{0,*}(\EDp*)$ are soft one
could speculate about the deformations
associated with $E$ as about the
Maurer-Cartan elements of the DGLA
$\OmEF^{0,*}(M,\EDp*)[[\h]]$\,.
Using the correspondence between
the Dolbeault and \v Cech pictures
one could relate these speculations to
Kontsevich's algebroid picture of
deformation quantization of
algebraic varieties \cite{K2}\,.

Finally, we think that the technique of
mixed resolutions proposed by A. Yekutieli
\cite{Ye} could help us to prove
Tsygan's formality conjecture for
Hochschild chains of the structure
sheaf of a smooth algebraic varieties
over an arbitrary field of
characteristic $0$\,.

\mbox{}

\noindent\footnotesize{\textsc{IRMA, 7 rue Ren\'e Descartes, F-67084 Strasbourg, France} \\
\emph{E-mail address}: {\bf calaque@math.u-strasbg.fr} \\
\emph{E-mail address}: {\bf halbout@math.u-strasbg.fr} \\

\noindent\textsc{Department of Mathematics, MIT, Cambridge, MA 02139, USA} \\
\emph{E-mail address}: {\bf vald@math.mit.edu}}

\end{document}